\documentclass[a4paper,12pt]{article}
\usepackage[all]{xy}
\usepackage[T1]{fontenc}
\usepackage[dvips]{graphicx}
\usepackage{amsmath,amssymb,amsfonts,amsthm,latexsym,mathrsfs,textcomp,verbatim}
\xyoption{web}

\title{Universal models and definability}
\author{{Olivia Caramello} \vspace{3 mm}\\ {\small DPMMS, University of Cambridge,}\\{\small Wilberforce Road, Cambridge CB3 0WB, UK}\\{\small O.Caramello@dpmms.cam.ac.uk}}
\date{June 17, 2009}
\begin{document}
\bgroup           
\let\footnoterule\relax  
\maketitle
\flushleft  
\begin{abstract}
We show that the investigation of universal models in Topos Theory can shed light on problems of definability in Logic as well as on the investigation of De Morgan's law and the law of excluded middle on Grothendieck toposes.     
\end{abstract} 
\egroup 
\flushleft
\vspace{5 mm}


\def\Monthnameof#1{\ifcase#1\or
   January\or February\or March\or April\or May\or June\or
   July\or August\or September\or October\or November\or December\fi}
\def\today{\number\day~\Monthnameof\month~\number\year}

%
%
%
\def\pushright#1{{
   \parfillskip=0pt            
   \widowpenalty=10000         
   \displaywidowpenalty=10000  
   \finalhyphendemerits=0      
  %
   \leavevmode                 
   \unskip                     
   \nobreak                    
   \hfil                       
   \penalty50                  
   \hskip.2em                  
   \null                       
   \hfill                      
   {#1}                        
  %
   \par}}                      

\def\qed{\pushright{$\square$}\penalty-700 \smallskip}

\newtheorem{theorem}{Theorem}[section]

\newtheorem{proposition}[theorem]{Proposition}

\newtheorem{scholium}[theorem]{Scholium}

\newtheorem{lemma}[theorem]{Lemma}

\newtheorem{corollary}[theorem]{Corollary}

\newtheorem{conjecture}[theorem]{Conjecture}

\newenvironment{proofs}%
 {\begin{trivlist}\item[]{\bf Proof }}%
 {\qed\end{trivlist}}

  \newtheorem{rmk}[theorem]{Remark}
\newenvironment{remark}{\begin{rmk}\em}{\end{rmk}}

  \newtheorem{rmks}[theorem]{Remarks}
\newenvironment{remarks}{\begin{rmks}\em}{\end{rmks}}

  \newtheorem{defn}[theorem]{Definition}
\newenvironment{definition}{\begin{defn}\em}{\end{defn}}

  \newtheorem{eg}[theorem]{Example}
\newenvironment{example}{\begin{eg}\em}{\end{eg}}

  \newtheorem{egs}[theorem]{Examples}
\newenvironment{examples}{\begin{egs}\em}{\end{egs}}


\mathcode`\<="4268  
\mathcode`\>="5269  
\mathcode`\.="313A  
\mathchardef\semicolon="603B 
\mathchardef\gt="313E
\mathchardef\lt="313C

\newcommand{\app}
 {{\sf app}}

\newcommand{\Ass}
 {{\bf Ass}}

\newcommand{\ASS}
 {{\mathbb A}{\sf ss}}

\newcommand{\Bb}
{\mathbb}

\newcommand{\biimp}
 {\!\Leftrightarrow\!}

\newcommand{\bim}
 {\rightarrowtail\kern-1em\twoheadrightarrow}

\newcommand{\bjg}
 {\mathrel{{\dashv}\,{\vdash}}}

\newcommand{\bstp}[3]
 {\mbox{$#1\! : #2 \bim #3$}}

\newcommand{\cat}
 {\!\mbox{\t{\ }}}

\newcommand{\cinf}
 {C^{\infty}}

\newcommand{\cinfrg}
 {\cinf\hy{\bf Rng}}

\newcommand{\cocomma}[2]
 {\mbox{$(#1\!\uparrow\!#2)$}}

\newcommand{\cod}
 {{\rm cod}}

\newcommand{\comma}[2]
 {\mbox{$(#1\!\downarrow\!#2)$}}

\newcommand{\comp}
 {\circ}

\newcommand{\cons}
 {{\sf cons}}

\newcommand{\Cont}
 {{\bf Cont}}

\newcommand{\ContE}
 {{\bf Cont}_{\cal E}}

\newcommand{\ContS}
 {{\bf Cont}_{\cal S}}

\newcommand{\cover}
 {-\!\!\triangleright\,}

\newcommand{\cstp}[3]
 {\mbox{$#1\! : #2 \cover #3$}}

\newcommand{\Dec}
 {{\rm Dec}}

\newcommand{\DEC}
 {{\mathbb D}{\sf ec}}

\newcommand{\den}[1]
 {[\![#1]\!]}

\newcommand{\Desc}
 {{\bf Desc}}

\newcommand{\dom}
 {{\rm dom}}

\newcommand{\Eff}
 {{\bf Eff}}

\newcommand{\EFF}
 {{\mathbb E}{\sf ff}}

\newcommand{\empstg}
 {[\,]}

\newcommand{\epi}
 {\twoheadrightarrow}

\newcommand{\estp}[3]
 {\mbox{$#1 \! : #2 \epi #3$}}

\newcommand{\ev}
 {{\rm ev}}

\newcommand{\Ext}
 {{\rm Ext}}

\newcommand{\fr}
 {\sf}

\newcommand{\fst}
 {{\sf fst}}

\newcommand{\fun}[2]
 {\mbox{$[#1\!\to\!#2]$}}

\newcommand{\funs}[2]
 {[#1\!\to\!#2]}

\newcommand{\Gl}
 {{\bf Gl}}

\newcommand{\hash}
 {\,\#\,}

\newcommand{\hy}
 {\mbox{-}}

\newcommand{\im}
 {{\rm im}}

\newcommand{\imp}
 {\!\Rightarrow\!}

\newcommand{\Ind}[1]
 {{\rm Ind}\hy #1}

\newcommand{\iten}[1]
{\item[{\rm (#1)}]}

\newcommand{\iter}
 {{\sf iter}}

\newcommand{\Kalg}
 {K\hy{\bf Alg}}

\newcommand{\llim}
 {{\mbox{$\lower.95ex\hbox{{\rm lim}}$}\atop{\scriptstyle
{\leftarrow}}}{}}

\newcommand{\llimd}
 {\lower0.37ex\hbox{$\pile{\lim \\ {\scriptstyle
\leftarrow}}$}{}}

\newcommand{\Mf}
 {{\bf Mf}}

\newcommand{\Mod}
 {{\bf Mod}}

\newcommand{\MOD}
{{\mathbb M}{\sf od}}

\newcommand{\mono}
 {\rightarrowtail}

\newcommand{\mor}
 {{\rm mor}}

\newcommand{\mstp}[3]
 {\mbox{$#1\! : #2 \mono #3$}}

\newcommand{\Mu}
 {{\rm M}}

\newcommand{\name}[1]
 {\mbox{$\ulcorner #1 \urcorner$}}

\newcommand{\names}[1]
 {\mbox{$\ulcorner$} #1 \mbox{$\urcorner$}}

\newcommand{\nml}
 {\triangleleft}

\newcommand{\ob}
 {{\rm ob}}

\newcommand{\op}
 {^{\rm op}}

\newcommand{\pepi}
 {\rightharpoondown\kern-0.9em\rightharpoondown}

\newcommand{\pmap}
 {\rightharpoondown}

\newcommand{\Pos}
 {{\bf Pos}}

\newcommand{\prarr}
 {\rightrightarrows}

\newcommand{\princfil}[1]
 {\mbox{$\uparrow\!(#1)$}}

\newcommand{\princid}[1]
 {\mbox{$\downarrow\!(#1)$}}

\newcommand{\prstp}[3]
 {\mbox{$#1\! : #2 \prarr #3$}}

\newcommand{\pstp}[3]
 {\mbox{$#1\! : #2 \pmap #3$}}

\newcommand{\relarr}
 {\looparrowright}

\newcommand{\rlim}
 {{\mbox{$\lower.95ex\hbox{{\rm lim}}$}\atop{\scriptstyle
{\rightarrow}}}{}}

\newcommand{\rlimd}
 {\lower0.37ex\hbox{$\pile{\lim \\ {\scriptstyle
\rightarrow}}$}{}}

\newcommand{\rstp}[3]
 {\mbox{$#1\! : #2 \relarr #3$}}

\newcommand{\scn}
 {{\bf scn}}

\newcommand{\scnS}
 {{\bf scn}_{\cal S}}

\newcommand{\semid}
 {\rtimes}

\newcommand{\Sep}
 {{\bf Sep}}

\newcommand{\sep}
 {{\bf sep}}

\newcommand{\Set}
 {{\bf Set }}

\newcommand{\Sh}
 {{\bf Sh}}

\newcommand{\ShE}
 {{\bf Sh}_{\cal E}}

\newcommand{\ShS}
 {{\bf Sh}_{\cal S}}

\newcommand{\sh}
 {{\bf sh}}

\newcommand{\Simp}
 {{\bf \Delta}}

\newcommand{\snd}
 {{\sf snd}}

\newcommand{\stg}[1]
 {\vec{#1}}

\newcommand{\stp}[3]
 {\mbox{$#1\! : #2 \to #3$}}

\newcommand{\Sub}
 {{\rm Sub}}

\newcommand{\SUB}
 {{\mathbb S}{\sf ub}}

\newcommand{\tbel}
 {\prec\!\prec}

\newcommand{\tic}[2]
 {\mbox{$#1\!.\!#2$}}

\newcommand{\tp}
 {\!:}

\newcommand{\tps}
 {:}

\newcommand{\tsub}
 {\pile{\lower0.5ex\hbox{.} \\ -}}

\newcommand{\wavy}
 {\leadsto}

\newcommand{\wavydown}
 {\,{\mbox{\raise.2ex\hbox{\hbox{$\wr$}
\kern-.73em{\lower.5ex\hbox{$\scriptstyle{\vee}$}}}}}\,}

\newcommand{\wbel}
 {\lt\!\lt}

\newcommand{\wstp}[3]
 {\mbox{$#1\!: #2 \wavy #3$}}

\newcommand{\fu}[2]
{[#1,#2]}

\newcommand{\st}[2]
 {\mbox{$#1 \to #2$}}

\tableofcontents

\newpage

\section{Introduction}
This paper is devoted to a general study of universal models in Topos Theory with a particular emphasis on their applications to definability by geometric formulae and to the investigation of the law of excluded middle and De Morgan's law on Grothendieck toposes.\\
In the second section, we establish a logical completeness result for universal models. This fact has many ramifications, some of which are explored in the third section. Here, among other results, we derive an explicit description of universal models of theories of presheaf type and for any such theory $\mathbb T$, we establish a duality between the category of finitely presented models of $\mathbb T$ and the full subcategory of the syntactic category of $\mathbb T$ on the formulae which present a $\mathbb T$-model; in this section, we also prove a definability theorem for theories of presheaf type.\\
The fourth section is devoted to describing the universal models for quotients ${\mathbb T}'$ of a theory of presheaf type $\mathbb T$. This is done for a quotient ${\mathbb T}'$ of $\mathbb T$ in terms of the Grothendieck topology on the opposite of the category of finitely presentable $\mathbb T$-models which corresponds to ${\mathbb T}'$ via the duality of Theorem 3.6 \cite{OC6}. These Grothendieck topologies are then described in several cases of interest, and accordingly natural descriptions of the universal models of the corresponding quotients are achieved; in this context, we also derive a sufficient condition for a quotient ${\mathbb T}'$ of $\mathbb T$ to have enough models.\\
In the fifth section, we discuss the relevance of universal models for the investigatation of definability of certain properties by geometric formulae; we observe that if $U$ is a universal model of a theory $\mathbb T$ over a signature $\Sigma$ in a topos $\cal E$ then every subobject of $UA_{1}\times \cdots \times UA_{n}$ (for any sorts $A_{1}, \ldots, A_{n}$ of $\Sigma$) is definable by a geometric formula and hence the lattice of subobjects of $UA_{1}\times \cdots \times UA_{n}$ in $\cal E$ is isomorphic to the lattice of subobjects of $\{x^{A_{1}}, \ldots, x^{A_{n}}. \top\}$ in the syntactic category ${\cal C}_{\mathbb T}$ of $\mathbb T$; as a consequence, the operation of pseudocomplementation in ${\cal C}_{\mathbb T}$, which plays a key role in our definability issues, acquires a natural semantic interpretation in terms of the universal model of $\mathbb T$. This leads us to establish a connection between definability of certain properties and validity of certain formulae in the universal model. The validity of such formulae is then shown to be related to the law of excluded middle and to De Morgan's law on toposes, and a thorough analysis of all these connections is carried out in the sixth section of the paper.\\  
In the last section, we present some applications in Algebra of the theory developed in the paper.

\section{Logical completeness of universal models}

Concerning notation, given two Grothendieck toposes $\cal E$ and $\cal F$ and a Grothendieck topology $J$ on a small category $\cal C$, we denote by ${\bf Geom}({\cal E}, {\cal F})$ the category of geometric morphisms from $\cal E$ to $\cal F$ and by ${\bf Flat}({\cal C}, {\cal E})$ (resp. ${\bf Flat}_{J}({\cal C}, {\cal E})$) the category of flat functors (resp. $J$-continuous flat functors) from $\cal C$ to $\cal E$. We denote by $\mathfrak{BTop}$ the $2$-category of Grothendieck toposes, geometric morphisms and geometric transformations between them.

\begin{definition}
Let $\mathbb T$ be a geometric theory. A $\mathbb T$-model $U$ in a Grothendieck topos $\cal G$ is said to be a universal model of $\mathbb T$ if for any $\mathbb T$-model $M$ in a Grothendieck topos $\cal F$ there exists a unique (up to isomorphism) geometric morphism $f:{\cal F}\to {\cal G}$ such that $f^{\ast}(U)\cong M$.
\end{definition}

\begin{rmk}\label{transferuniv}
\emph{We note that, by the $2$-dimensional Yoneda Lemma, if a topos $\cal G$ contains a universal model of a geometric theory $\mathbb T$ then $\cal G$ satisfies the universal property of the classifying topos of $\mathbb T$. Conversely, if a topos $\cal E$ classifies a geometric theory $\mathbb T$ then $\cal E$ contains a universal model of $\mathbb T$. It also follows that classifying toposes, and hence universal models, are uniquely determined up to equivalence. More precisely, if $M$ and $N$ are universal models of a geometric theory $\mathbb T$ lying respectively in toposes $\cal F$ and $\cal G$ then there exists a unique (up to isomorphism) geometric equivalence between $\cal F$ and $\cal G$ such that its inverse image functors send $M$ and $N$ to each other (up to isomorphism). Indeed, by definition of universal model, there are equivalences $\tau^{\cal E}_{M}:{\bf Geom}({\cal E}, {\cal F}) \simeq {\mathbb T}\textrm{-mod}({\cal E})$ and $\tau^{\cal E}_{N}:{\bf Geom}({\cal E}, {\cal G}) \simeq {\mathbb T}\textrm{-mod}({\cal E})$ natural in ${\cal E}\in \mathfrak{Btop}$ such that $\tau_{M}^{\cal F}(1_{\cal F})=M$ and $\tau_{N}^{\cal G}(1_{\cal G})=N$. Then by composing the equivalences $\tau_{M}$ and $\tau_{N}$ we get, by the 2-dimensional Yoneda Lemma, a unique up to isomorphism geometric equivalence between $\cal F$ and $\cal G$ such that its inverse image functors sends $M$ and $N$ to each other (up to isomorphism), as required.}
\end{rmk}
The following fact will be exploited extensively in this paper.

\begin{lemma}\label{subtopos}
Let $\cal E$ be the classifying topos of a geometric theory $\mathbb T$ and $i:{\cal F}\hookrightarrow {\cal E}$ a subtopos of $\cal E$ with associated sheaf functor $a_{{\cal F}}:{\cal E}\to {\cal F}$ which corresponds to a quotient ${\mathbb T}'$ of $\mathbb T$ via Theorem 3.6 \cite{OC6}. If $U$ is a universal model of $\mathbb T$ in $\cal E$ then $a_{{\cal F}}(U)$ is a universal model of ${\mathbb T}'$ in $\cal F$. 
\end{lemma}

\begin{proofs}
The thesis follows at once from the commutativity and naturality in ${\cal G}\in \mathfrak{BTop}$ of the diagram
\[  
\xymatrix {
{{\mathbb T}'}\emph{\textrm{-mod}}({\cal G})  \ar[d]^{i^{\cal G}_{{\mathbb T}'}} & & \ar[ll]^{\simeq} {\bf Geom}({\cal G}, {\cal F}) \ar[d]^{i \circ -} \\
{\mathbb T}\emph{\textrm{-mod}}({\cal G})  & &  \ar[ll]^{\simeq}_{\tau^{\cal G}} {\bf Geom}({\cal G}, {\cal E})}
\]
where $i^{\cal G}_{{\mathbb T}'}$ is the obvious inclusion and the equivalence $\tau^{\cal G}$ sends a geometric morphism $f$ in ${\bf Geom}({\cal G}, {\cal E})$ to the $\mathbb T$-model $f^{\ast}(U)$ (see Remark 3.7 \cite{OC6}).  
\end{proofs}

Let us recall from section D3.1 \cite{El2} that, given a geometric theory $\mathbb T$ over a signature $\Sigma$, we have a syntactic description of a universal model $U_{\mathbb T}$ of $\mathbb T$ lying in the classifying topos $\Sh({\cal C}_{\mathbb T}, J_{\mathbb T})$ of $\mathbb T$, where ${\cal C}_{\mathbb T}$ is the geometric syntactic category of $\mathbb T$ and $J_{\mathbb T}$ is the canonical topology on ${\cal C}_{\mathbb T}$. Specifically, $U_{\mathbb T}$ is the image of the `universal' model $M_{\mathbb T}$ of $\mathbb T$ in ${\cal C}_{\mathbb T}$ under the Yoneda embedding $y:{\cal C}_{\mathbb T}\hookrightarrow \Sh({\cal C}_{\mathbb T}, J_{\mathbb T})$; recall that $M_{\mathbb T}$ assigns to a sort $A$ the object $\{x^{A}.\top\}$ where $x^{A}$ is a variable of sort $A$, to a function symbol $f:A_{1}\cdots A_{n}\to B$ the morphism
\[  
\xymatrix {
\{x_{1}^{A_{1}}, \ldots , x_{n}^{A_{n}}. \top\} \ar[rrrr]^{[f(x_{1}^{A_{1}}, \ldots , x_{n}^{A_{n}})=y^{B}]} & & & & \{y^{B}.\top\} }
\]
and to a relation symbol $R\mono A_{1}\cdots A_{n}$ the subobject 
\[  
\xymatrix {
\{x_{1}^{A_{1}}, \ldots , x_{n}^{A_{n}}. R(x_{1}^{A_{1}}, \ldots , x_{n}^{A_{n}})\} \ar[rrrr]^{[R(x_{1}^{A_{1}}, \ldots , x_{n}^{A_{n}})]} & & & & \{x_{1}^{A_{1}}, \ldots , x_{n}^{A_{n}}. \top\} }
\]
as in Lemma D1.4.4(iv) \cite{El2}.\\
Moreover, we recall that $U_{\mathbb T}$ is a conservative model of $\mathbb T$ (i.e. for any geometric sequent $\sigma$ over $\Sigma$, if $\sigma$ is valid in $U_{\mathbb T}$ then $\sigma$ is provable in $\mathbb T$ using geometric logic), and hence, by Remark \ref{transferuniv}, any universal model $M$ of a geometric theory $\mathbb T$ is a conservative model of $\mathbb T$.\\

Universal models of geometric theories enjoy a form of logical completeness; specifically, we have the following result.
\begin{theorem}\label{definab}
Let $\mathbb T$ be a geometric theory over a signature $\Sigma$ and $U$ a universal model of $\mathbb T$ in a topos $\cal E$. Then:\\
(i) for any subobject $S\mono UA_{1}\times \cdots UA_{n}$ in $\cal E$ there exists a geometric formula $\phi(x^{A_{1}}, \ldots, x^{A_{n}})$ over $\Sigma$ such that $S=[[\phi]]_{U}$;\\
(ii) for any arrow $f:[[\vec{x}.\phi]]_{U}\to [[\vec{y}.\psi]]_{U}$ in $\cal E$ where $\phi(\vec{x})$ and $\psi(\vec{y})$ are geometric formulae over $\Sigma$, there exists a geometric formula $\theta(\vec{x}, \vec{y})$ over $\Sigma$ such that the sequents $(\phi \vdash_{\vec{x}} (\exists y)\theta)$, $(\theta \vdash_{\vec{x}, \vec{y}} \phi \wedge \psi)$ and $((\theta \wedge \theta[\vec{y'}\slash \vec{y}]) \vdash_{\vec{x}, \vec{y}, \vec{y'}} (\vec{y}=\vec{y'}))$ are provable in $\mathbb T$ and $[[\theta]]_{U}$ is the graph of $f$. 
\end{theorem}
 
\begin{proofs}
First, we note that, by Remark \ref{transferuniv} and the fact that inverse image functors of geometric morphisms preserve the interpretation of geometric formulae, we can suppose without loss of generality that $U$ is the universal model $U_{\mathbb T}$ of $\mathbb T$ lying in the classifying topos $\Sh({\cal C}_{\mathbb T}, J_{\mathbb T})$ of $\mathbb T$ described above.\\
To prove (i), we observe that, given a geometric formula $\phi(x^{A_{1}}, \ldots, x^{A_{n}})$ over $\Sigma$, the interpretation $[[\phi]]_{U_{\mathbb T}}$ of $\phi$ in $U_{\mathbb T}$ identifies with the ($J_{\mathbb T}$-closed) sieve on $\{x^{A_{1}}, \ldots, x^{A_{n}}. \top\}$ generated by the canonical monomorphism $[\phi]:\{\vec{x}. \phi\}\mono \{x^{A_{1}}, \ldots, x^{A_{n}}. \top\}$. Now, if $S$ is a subobject of $U_{\mathbb T}A_{1}\times \cdots U_{\mathbb T}A_{n}\cong Hom_{{\cal C}_{\mathbb T}}(-, \{x^{A_{1}}, \ldots, x^{A_{n}}. \top\})$ then $S$, regarded as a sieve, is $J_{\mathbb T}$-closed and hence, by Proposition 2.6 \cite{OC6}, it is generated by a subobject of $\{x^{A_{1}}, \ldots, x^{A_{n}}. \top\}$ in ${\cal C}_{\mathbb T}$. Thus the thesis follows from the characterization of subobjects in ${\cal C}_{\mathbb T}$ given by Lemma D1.4.4(iv) \cite{El2}.\\
Let us now prove (ii). By the Yoneda Lemma, any arrow $f:[[\vec{x}. \phi]]_{{\mathbb U}_{\mathbb T}}\cong Hom_{{\cal C}_{\mathbb T}}(-, \{\vec{x}. \phi\})\to Hom_{{\cal C}_{\mathbb T}}(-, \{\vec{y}. \psi\})\cong [[\vec{y}.\psi]]_{{\mathbb U}_{\mathbb T}}$ in $\Sh({\cal C}_{\mathbb T}, J_{\mathbb T})$ is of the form $Hom_{{\cal C}_{\mathbb T}}(-,[\theta])$ for a unique arrow $[\theta]:\{\vec{x}. \phi\}\to \{\vec{y}.\psi\}$ in ${\cal C}_{\mathbb T}$; but $[[\theta]]_{U_{\mathbb T}}$ is the graph of the arrow $Hom_{{\cal C}_{\mathbb T}}(-,[\theta])$, from which our thesis follows.     
\end{proofs} 
\begin{rmk}\label{restr}
\emph{It is natural to wonder how the theorem can be adapted to smaller fragments of logic. It is clear that if $\mathbb T$ is a cartesian (resp. regular, coherent) theory then part (ii) holds with geometric replaced by cartesian (resp. regular, coherent) everywhere (use the appropriate syntactic representation of the classifying topos of $\mathbb T$). However, part (i) no longer survives.}
\end{rmk}

The following proposition gives a useful property of universal models.

\begin{proposition}\label{separ}
Let $\mathbb T$ be a geometric theory over a signature $\Sigma$ and $U$ a universal model of $\mathbb T$ in a Grothendieck topos $\cal E$. Then the subobjects in $\cal E$ of the objects of the form $UA_{1}\times \cdots \times UA_{n}$ for sorts $A_{1}, \ldots, A_{n}$ of $\Sigma$ form a separating set of $\cal E$.
\end{proposition}

\begin{proofs}
By using Remark \ref{transferuniv}, it suffices to prove the thesis for $U$ equal to the universal model $U_{\mathbb T}$ of $\mathbb T$ in the topos $\Sh({\cal C}_{\mathbb T}, J_{\mathbb T})$. Now, every geometric formula $\phi(x^{A_{1}}, \ldots, x^{A_{n}})$ gives rise to a subobject $[[\phi]]_{U}$  of $UA_{1}\times \cdots \times UA_{n}$ (the converse is also true by Theorem \ref{definab}). But if $U$ equal to the universal model $U_{\mathbb T}$ then $[[\phi]]_{U}$ is given by $y(\{\vec{x}.\phi\})$ where $y:{\cal C}_{\mathbb T}\to \Sh({\cal C}_{\mathbb T}, J_{\mathbb T})$ is the Yoneda embedding, and hence our thesis follows from the fact that the objects of the form $y(\{\vec{x}.\phi\})$ form a separating set for $\Sh({\cal C}_{\mathbb T}, J_{\mathbb T})$.    
\end{proofs}

\newpage
\section{Universal models of theories of presheaf type}\label{sec}
We recall that a theory of presheaf type is a geometric theory classified by a presheaf topos (equivalently, by the topos $[\textrm{f.p.} {\mathbb T}\textrm{-mod}(\Set), \Set]$, where $\textrm{f.p.} {\mathbb T}\textrm{-mod}(\Set)$ is the full subcategory of ${\mathbb T}\textrm{-mod}(\Set)$ on the finitely presentable objects). In particular, every cartesian theory is of presheaf type (cfr. Corollary D3.1.2 \cite{El2}).\\
Let $\mathbb T$ be a theory of presheaf type, together with an equivalence $\xi^{\cal E}: {\bf Flat}(\textrm{f.p.} {\mathbb T}\textrm{-mod}(\Set)^{\textrm{op}}, {\cal E}) \to {\mathbb T}\textrm{-mod}({\cal E})$ natural in ${\cal E}\in \mathfrak{Btop}$. If $y:\textrm{f.p.} {\mathbb T}\textrm{-mod}(\Set) \to {\bf Flat}(\textrm{f.p.} {\mathbb T}\textrm{-mod}(\Set)^{\textrm{op}}, \Set)$ is the Yoneda embedding then the factorization of the composite $\xi^{\Set}\circ y:\textrm{f.p.} {\mathbb T}\textrm{-mod}(\Set) \to {\mathbb T}\textrm{-mod}({\Set})$ through the inclusion $i:\textrm{f.p.} {\mathbb T}\textrm{-mod}(\Set)\hookrightarrow {\mathbb T}\textrm{-mod}(\Set)$ is an equivalence of categories $\tau^{\xi}:\textrm{f.p.} {\mathbb T}\textrm{-mod}(\Set) \to \textrm{f.p.} {\mathbb T}\textrm{-mod}(\Set)$. The equivalence $\xi$ is said to be \emph{canonical} if $\tau^{\xi}(c)\cong c$ naturally in $c\in \textrm{f.p.} {\mathbb T}\textrm{-mod}(\Set)$. We observed in \cite{OC6} that, given a theory of presheaf type $\mathbb T$, there is always a canonical equivalence $\xi^{\cal E}: {\bf Flat}(\textrm{f.p.} {\mathbb T}\textrm{-mod}(\Set)^{\textrm{op}}, {\cal E}) \to {\mathbb T}\textrm{-mod}({\cal E})$ natural in ${\cal E}\in \mathfrak{Btop}$.\\  
The following result gives an explicit description of `the' universal model of a given theory of presheaf type. 

\begin{theorem}\label{presheaf}
Let $\mathbb T$ be a theory of presheaf type over a signature $\Sigma$. Then the $\Sigma$-structure $M_{\mathbb T}$ in $[\textrm{f.p.} {\mathbb T}\textrm{-mod}(\Set), \Set]$ which assigns to a sort $A$ the functor $M_{\mathbb T}A$ given by $(M_{\mathbb T}A)(M)=MA$, to a function symbol $f:A_{1}\cdots A_{n}\to B$ the morphism $M_{\mathbb T}A_{1}\times \cdots M_{\mathbb T}A_{n}\to M_{\mathbb T}B$ given by $(M_{\mathbb T}f)(M)=Mf$ and to a relation symbol $R\mono A_{1}\cdots A_{n}$ the subobject $M_{\mathbb T}R\mono M_{\mathbb T}A_{1}\times \cdots M_{\mathbb T}A_{n}$ given by $(M_{\mathbb T}R)(M)=MR$ (for any $M\in \textrm{f.p.} {\mathbb T}\textrm{-mod}(\Set)$) is a universal model for $\mathbb T$; moreover, for any geometric formula $\phi(\vec{x})$ over $\Sigma$, the interpretation $[[\vec{x}. \phi]]_{M_{\mathbb T}}$ of $\phi(\vec{x})$ in $M_{\mathbb T}$ is given by $[[\vec{x}. \phi]]_{M_{\mathbb T}}(M)=[[\vec{x}. \phi]]_{M}$ for any $M\in \textrm{f.p.} {\mathbb T}\textrm{-mod}(\Set)$. In particular, the finitely presentable $\mathbb T$-models are jointly conservative for $\mathbb T$.    
\end{theorem}
\begin{proofs}
Given a small category $\cal C$ and an object $c\in {\cal C}$, since in functor categories small limits and colimits are computed pointwise, the evaluation functor $ev_{c}:[{\cal C}, \Set]\to \Set$ at $c$ preserves all small limits and colimits and hence, by Remark C2.2.10 \cite{El2}, it is the inverse image functor of a geometric morphism $e_{c}:\Set \to [{\cal C}, \Set]$. By using the adjunction $(e_{c})_{\ast} \dashv ev_{c}$ and the Yoneda Lemma, one gets the following explicit description of the direct image functor $(e_{c})_{\ast}:\Set \to [{\cal C}, \Set]$: $(e_{c})_{\ast}(A)(c')=Hom_{\Set}({\cal C}(c,c'), A)$ for any set $A$ and object $c'\in {\cal C}$.\\
Now, let us consider a canonical equivalence $\xi^{\cal E}: {\bf Flat}(\textrm{f.p.} {\mathbb T}\textrm{-mod}(\Set)^{\textrm{op}}, {\cal E}) \to {\mathbb T}\textrm{-mod}({\cal E})$ natural in ${\cal E}\in \mathfrak{Btop}$ for the theory $\mathbb T$. By composing it with Diaconescu's equivalence ${\bf Geom}({\cal E}, [\textrm{f.p.} {\mathbb T}\textrm{-mod}(\Set), \Set]) \simeq {\bf Flat}(\textrm{f.p.} {\mathbb T}\textrm{-mod}(\Set)^{\textrm{op}}, {\cal E})$, we get an equivalence $\tau^{\cal E}:{\bf Geom}({\cal E}, [\textrm{f.p.} {\mathbb T}\textrm{-mod}(\Set), \Set]) \simeq {\mathbb T}\textrm{-mod}({\cal E})$ natural in ${\cal E}\in \mathfrak{Btop}$; let us define $U$ to be the image of the identical geometric morphism on $[\textrm{f.p.} {\mathbb T}\textrm{-mod}(\Set), \Set]$ under this equivalence. Then $U$ is a universal model of $\mathbb T$ and a geometric morphism $f:{\cal E}\to [\textrm{f.p.} {\mathbb T}\textrm{-mod}(\Set), \Set]$ is sent under $\tau^{\cal E}$ to the $\mathbb T$-model $f^{\ast}(U)$. Recall that, via Diaconescu's equivalence, the identity on $[\textrm{f.p.} {\mathbb T}\textrm{-mod}(\Set), \Set]$ is sent to the flat functor $y:\textrm{f.p.} {\mathbb T}\textrm{-mod}(\Set)^{\textrm{op}} \to [\textrm{f.p.} {\mathbb T}\textrm{-mod}(\Set), \Set]$ given by the Yoneda embedding and hence, by naturality, the geometric morphism $e_{M}:\Set \to [\textrm{f.p.} {\mathbb T}\textrm{-mod}(\Set), \Set]$ (for $M\in \textrm{f.p.} {\mathbb T}\textrm{-mod}(\Set)$) is sent to the flat functor $(e_{M})^{\ast}\circ y=Hom_{\textrm{f.p.} {\mathbb T}\textrm{-mod}(\Set)}(-, M): \textrm{f.p.} {\mathbb T}\textrm{-mod}(\Set)^{\textrm{op}} \to \Set$. But, since $\xi$ is canonical, $Hom_{\textrm{f.p.} {\mathbb T}\textrm{-mod}(\Set)}(-, M)$ is sent to $M$ via $\xi^{\Set}$, and hence $\tau^{\Set}$ sends $e_{M}$ to $M$ (for any $M\in \textrm{f.p.} {\mathbb T}\textrm{-mod}(\Set)$). Thus $M\cong (e_{M})^{\ast}(U)=ev_{M}(U)$ for any $M\in \textrm{f.p.} {\mathbb T}\textrm{-mod}(\Set)$ and hence $U$ is (isomorphic to) the $\Sigma$-structure $M_{\mathbb T}$ defined in the statement of the theorem.\\
Now, the fact that $[[\vec{x}. \phi]]_{M_{\mathbb T}}(M)=[[\vec{x}. \phi]]_{M}$ for any geometric formula $\phi(\vec{x})$ and model $M\in \textrm{f.p.} {\mathbb T}\textrm{-mod}(\Set)$ follows from the fact that the functors $ev_{M}$ are geometric (being inverse image functors of geometric morphisms).   \end{proofs}
We note that our theorem specializes to Corollary D3.1.2 \cite{El2} in the case $\mathbb T$ is cartesian.\\
It is natural to wonder what the interpretations of arbitrary first-order formulae in the universal model $M_{\mathbb T}$ look like; in particular, one can ask if we also have $[[\vec{x}. \phi]]_{M_{\mathbb T}}(M)=[[\vec{x}. \phi]]_{M}$ for any first-order formula $\phi(\vec{x})$ and model $M\in \textrm{f.p.} {\mathbb T}\textrm{-mod}(\Set)$. This condition is satisfied if the evaluation functors $ev_{M}$ are Heyting functors i.e. the geometric morphisms $e_{M}$ are open; but it is immediate to see, by using Lemma C3.1.2 \cite{El2}, that this is the case if and only if the category $\textrm{f.p.} {\mathbb T}\textrm{-mod}(\Set)$ is a groupoid.\\

The following result is a definability theorem for theories of presheaf type.

\begin{corollary}\label{predef}
Let $\mathbb T$ be a theory of presheaf type over a signature $\Sigma$, $A_{1}, \ldots, A_{n}$ a string of sorts of $\Sigma$ and suppose we are given, for every finitely presentable $\Set$-model $M$ of $\mathbb T$ a subset $R_{M}$ of $MA_{1}\times \ldots \times MA_{n}$ in such a way that each $\mathbb T$-model homomorphism $h:M\to N$ maps $R_{M}$ into $R_{N}$. Then there exists a geometric formula-in-context $\phi(x^{A_{1}}, \ldots, x^{A_{n}})$ such that $R_{M}=[[\phi]]_{M}$ for each $M$.       
\end{corollary}

\begin{proofs}
The corollary immediately follows from Theorem \ref{definab}(i) and Theorem \ref{presheaf}, observing that the assignment $M\to R_{M}$ in the statement of the corollary gives rise to a subfunctor $R\mono M_{\mathbb T}$, where $M_{\mathbb T}$ is the universal model of $\mathbb T$ in $[\textrm{f.p.} {\mathbb T}\textrm{-mod}(\Set), \Set]$ defined in the statement of Theorem \ref{presheaf}. 
\end{proofs}  
Note that, conversely, for any geometric formula $\phi$, any $\mathbb T$-model homomorphism $h:M\to N$ maps $[[\phi]]_{M}$ into $[[\phi]]_{N}$ (cfr. Lemma D1.2.9 \cite{El2}).\\ 

Let us now recall from \cite{OC6} the following notions.\\
\begin{definition}
Let $\mathbb T$ be a geometric theory over a signature $\Sigma$ and $\phi(x_{1}^{A_{1}}, \ldots, x_{n}^{A_{n}})$ a geometric formula over $\Sigma$. We say that a $\mathbb T$-model $M$ in \Set is finitely presented by $\phi$ (or that $\phi$ presents $M$) if there exists a string of elements $\vec{\xi_{\phi}}:=
(\xi_{1}, \ldots, \xi_{n})\in MA_{1}\times \ldots MA_{n}$, called the generators of $M$, such that for any $\mathbb T$-model $N$ in \Set and string of elements $(b_{1}, \ldots, b_{n})\in MA_{1}\times \ldots MA_{n}$ such that $(b_{1}, \ldots, b_{n})\in [[\phi]]_{N}$, there exists a unique arrow $f_{(b_{1}, \ldots, b_{n})}:M\to N$ in ${\mathbb T}\textrm{-mod}(\Set)$ such that $f(\vec{\xi_{\phi}}):=(fA_{1}\times \ldots fA_{n})((\xi_{1}, \ldots, \xi_{n}))=(b_{1}, \ldots, b_{n})$.      
\end{definition}
Of course, there can be at most one (up to isomorphism) $\mathbb T$-model finitely presented by a given formula-in-context $\{\vec{x}.\phi\}$; we will denote such model by $M_{\{\vec{x}.\phi\}}$, or simply $M_{\phi}$ when the context $\vec{x}$ can be obviously inferred.\\
Given a geometric theory $\mathbb T$ over a signature $\Sigma$, a geometric formula $\phi(x_{1}^{A_{1}}, \ldots, x_{n}^{A_{n}})$ over $\Sigma$ and a Grothendieck topos $\cal E$, consider the functor $F_{\phi}^{\cal E}: {\mathbb T}\textrm{-mod}({\cal E}) \to {\cal E}$ which sends to each model $N\in {\mathbb T}\textrm{-mod}({\cal E})$ (the domain of) the interpretation $[[\phi]]_{N}$ of $\phi$ in $N$ and acts on arrows in the obvious way. Note that, by the soundness theorem for geometric logic, for any topos $\cal E$ the assignment $\phi \to F^{\cal E}_{\phi}$ is functorial on the geometric syntactic category of $\mathbb T$ i.e. it gives rise to a functor $F_{\mathbb T}^{\cal E}:{\cal C}_{\mathbb T}\to [{\mathbb T}\textrm{-mod}({\cal E}), {\cal E}]$.\\
The functor $F^{\Set}_{\phi}$ preserves filtered colimits (cfr. the proof of Lemma D2.4.9 \cite{El2}) so if it is representable then the representing object is a finitely presentable model. Notice that, by the Yoneda Lemma, $F^{\Set}_{\phi}$ is representable if and only if there exists a $\mathbb T$-model finitely presented by $\phi$.\\ 
The functor $M_{\mathbb T}A: \textrm{f.p.} {\mathbb T}\textrm{-mod}(\Set) \to \Set$ considered above is the restriction to $\textrm{f.p.} {\mathbb T}\textrm{-mod}(\Set)$ of the functor $F^{\Set}_{{\top(x^{A})}}$; so if the formula $\{x^{A}. \top\}$ presents a $\mathbb T$-model $M_{\{x^{A}.\top\}}$ then $M_{\mathbb T}A$ is representable with representing object $M_{\{x^{A}.\top\}}$.\\
Given a theory of presheaf type $\mathbb T$ with geometric syntactic category ${\cal C}_{\mathbb T}$, let us denote by $\tilde{{\cal C}_{\mathbb T}}$ the full subcategory of ${\cal C}_{\mathbb T}$ on the formulae-in-context which present a $\mathbb T$-model, and by $\textrm{f.p.} {\mathbb T}\textrm{-mod}(\Set)_{s}$ the full subcategory of $\textrm{f.p.} {\mathbb T}\textrm{-mod}(\Set)$ on the finitely presented models.\\
We can define a `dualizing' functor $d:\tilde{{\cal C}_{\mathbb T}}\to \textrm{f.p.} {\mathbb T}\textrm{-mod}(\Set)^{\textrm{op}}$ with values in ${\textrm{f.p.} {\mathbb T}\textrm{-mod}(\Set)}_{s}^{\textrm{op}}$ as follows: given $\{\vec{x}. \phi\}\in \tilde{{\cal C}_{\mathbb T}}^{\textrm{op}}$, $d(\{\vec{x}. \phi\})=M_{\{\vec{x}.\phi\}}$ and given an arrow $[\theta]:\{\vec{y}. \psi\}\to \{\vec{x}. \phi\}$ in $\tilde{{\cal C}_{\mathbb T}}$, $d([\theta]):M_{\{\vec{x}.\phi\}} \to M_{\{\vec{y}.\psi\}}$ is the arrow in $\textrm{f.p.}{\mathbb T}\textrm{-mod}(\Set)$ defined by setting $d([\theta])(\vec{{\xi}_{\phi}})=[[\theta]]_{M_{\psi}}(\vec{\xi_{\psi}})$.\\

From Theorems \ref{definab} and \ref{presheaf}, we know that any sieve $R$ in $\textrm{f.p.} {\mathbb T}\textrm{-mod}(\Set)^{\textrm{op}}$ on a finitely presented $\mathbb T$-model $M_{\{\vec{x}.\phi\}}$ is, regarded as a subfunctor $R\mono Hom_{\textrm{f.p.} {\mathbb T}\textrm{-mod}(\Set)}(M_{\phi}, -)=[[\vec{x}. \phi]]_{M_{\mathbb T}}$, the interpretation in the universal model $M_{\mathbb T}$ of $\mathbb T$ of a geometric formula $\psi(\vec{x})$ over $\Sigma$ such that $\psi \vdash_{\vec{x}} \phi$ is provable in $\mathbb T$; moreover, $R$ identifies with the collection of arrows $f:M_{\phi}\to M$ in $\textrm{f.p.} {\mathbb T}\textrm{-mod}(\Set)$ such that $f(\vec{\xi_{\phi}})\in [[\vec{x}. \psi]]_{M}$. The following result provides an explicit such formula $\psi(\vec{x})$ for sieves $R$ of the form $d(S)$ for some sieve $S$ in $\tilde{{\cal C}_{\mathbb T}}$ on $\{\vec{x}. \phi\}$.  

\begin{theorem}\label{dualint}
Let $\mathbb T$ be a theory of presheaf type and $S$ a sieve in $\tilde{{\cal C}_{\mathbb T}}$ on $\{\vec{x}. \phi\}$ generated by arrows $[\theta_{i}]:\{\vec{y_{i}}. \psi\}\to \{\vec{x}. \phi\}$ as $i\in I$. Then, with the notation above, the sieve $d(S)$ in $\textrm{f.p.} {\mathbb T}\textrm{-mod}(\Set)$ generated by the arrows $d(f)$ as $f$ varies in $S$, regarded as a functor $\textrm{f.p.} {\mathbb T}\textrm{-mod}(\Set)\to \Set$, is equal to $[[\mathbin{\mathop{\textrm{\huge $\vee$}}\limits_{i\in I}}(\exists \vec{y_{i}})\theta_{i}]]_{M_{\mathbb T}}$.   
\end{theorem}

\begin{proofs}
Since for any sieve $S$ on an object $c$ in a category $\cal C$, the corresponding subobject $S\mono {\cal C}(-, c)$ in $[{\cal C}^{\textrm{op}}, \Set]$ is the union of the subobjects $(f)\mono {\cal C}(-, c)$ in $[{\cal C}^{\textrm{op}}, \Set]$ as $f$ varies in $\cal C$, it clearly suffices to prove the thesis for $S$ generated by a single arrow $[\theta]:\{\vec{y}. \psi\}\to \{\vec{x}. \phi\}$ in $\tilde{{\cal C}_{\mathbb T}}$. Hence, by Theorem \ref{presheaf}, we are reduced to prove that for any arrow $f:M_{\phi}\to M$ in ${\textrm{f.p.} {\mathbb T}\textrm{-mod}(\Set)}$ with domain $M_{\phi}$, $f$ factors through $d([\theta])$ if and only if $f(\vec{\xi_{\phi}})\in [[\exists \vec{y}\theta]]_{M}$.\\
First, note that $[[\exists \vec{y}\theta]]_{M}$ is the image in $\Set$ of the arrow $[[\theta]]_{M}:[[\vec{y}. \psi]]_{M}\to [[\vec{x}. \phi]]_{M}$; so $f(\vec{\xi_{\phi}})\in [[\exists \vec{y}\theta]]_{M}$ if and only if there exists an element $\vec{a}\in [[\vec{y}. \psi]]_{M}$ such that $[[\theta]]_{M}(\vec{a})=f(\vec{\xi_{\phi}})$. But, by definition of $M_{\phi}$, $\vec{a}$ identifies with an arrow $g:M_{\psi}\to M$ in $\textrm{f.p.} {\mathbb T}\textrm{-mod}(\Set)$ such that $g(\vec{\xi_{\psi}})=\vec{a}$ and hence, under this identification, the condition $f(\vec{\xi_{\phi}})\in [[\exists \vec{y}\theta]]_{M}$ rewrites as follows: there exists an arrow $g:M_{\psi}\to M$ in $\textrm{f.p.} {\mathbb T}\textrm{-mod}(\Set)$ such that $[[\theta]]_{M}(g(\vec{\xi_{\psi}}))=f(\vec{\xi_{\phi}})$. Now, by the functoriality of $F_{\mathbb T}^{[\textrm{f.p.} {\mathbb T}\textrm{-mod}(\Set), \Set]}$, we have $[[\theta]]_{M}(g(\vec{\xi_{\psi}}))=g([[\theta]]_{M}(\vec{\xi_{\psi}}))=g(d([\theta])(\vec{\xi_{\phi}}))=(g \circ d([\theta]))(\vec{\xi_{\phi}})$, from which our thesis follows.  
\end{proofs}
\begin{rmk}\label{projection}
\emph{As an instance of Theorem \ref{dualint}, we note that if $\{\vec{x}. \top\}$ and $\{\vec{x}. \phi\}$ are objects of $\tilde{{\cal C}_{\mathbb T}}$ then the cosieve in $\textrm{f.p.} {\mathbb T}\textrm{-mod}(\Set)$ on $M_{\{\vec{x}.\top\}}$ generated by the arrow $d([\phi])$ where $[\phi]$ is the canonical monomorphism $\{\vec{x}. \phi\}\mono \{\vec{x}. \top\}$ in ${\cal C}_{\mathbb T}$ identifies with $[[\phi]]_{M_{\mathbb T}}$.}
\end{rmk}
 
\begin{theorem}\label{dualitap}
With the notation above, for any theory of presheaf type $\mathbb T$ the functor $d:\tilde{{\cal C}_{\mathbb T}}\to {\textrm{f.p.} {\mathbb T}\textrm{-mod}(\Set)}_{s}^{\textrm{op}}$ is an equivalence of categories.
\end{theorem}

\begin{proofs}
Consider the composite $y \circ d:\tilde{{\cal C}_{\mathbb T}} \to [\textrm{f.p.} {\mathbb T}\textrm{-mod}(\Set), \Set]$ of $d:\tilde{{\cal C}_{\mathbb T}}\to {\textrm{f.p.} {\mathbb T}\textrm{-mod}(\Set)^{\textrm{op}}}$ with the Yoneda embedding $y:{\textrm{f.p.} {\mathbb T}\textrm{-mod}(\Set)^{\textrm{op}}} \to [\textrm{f.p.} {\mathbb T}\textrm{-mod}(\Set), \Set]$.\\
If $[\theta]:\{\vec{y}. \psi\}\to \{\vec{x}. \phi\}$ is an arrow in $\tilde{{\cal C}_{\mathbb T}}$ then $(y\circ d)([\theta])$ is equal to $[[\theta]]_{M_{\mathbb T}}:[[\vec{y}. \psi]]_{M_{\mathbb T}}\to [[\vec{x}. \phi]]_{M_{\mathbb T}}$. To prove this, by Theorem \ref{presheaf} it suffices to prove that for any $M\in {\textrm{f.p.} {\mathbb T}\textrm{-mod}(\Set)}$, the diagram
\[  
\xymatrix {
[[\vec{y}. \psi]]_{M} \ar[rr]^{[[\theta]]_{M}} \ar[d]^{\cong} & & \ar[d]^{\cong} [[\vec{x}. \phi]]_{M}\\
Hom_{{\textrm{f.p.} {\mathbb T}\textrm{-mod}(\Set)}}(M_{\psi}, M) \ar[rr]^{- \circ d([\theta])} & & Hom_{{\textrm{f.p.} {\mathbb T}\textrm{-mod}(\Set)}}(M_{\phi}, M) }
\]
commutes. But this is equivalent to saying that for any arrow $f:M_{\psi}\to M$ in $\textrm{f.p.} {\mathbb T}\textrm{-mod}(\Set)$, $[[\theta]]_{M}(f(\vec{\xi_{\psi}}))=(f\circ d([\theta]))(\vec{\xi_{\phi}})$, and this condition is satisfied since, by the functoriality of $F_{\mathbb T}^{[\textrm{f.p.} {\mathbb T}\textrm{-mod}(\Set), \Set]}$, $[[\theta]]_{M}(f(\vec{\xi_{\psi}}))=f([[\theta]]_{M_{\psi}}(\vec{\xi_{\psi}}))$.\\
Hence, $y$ being full and faithful, $d$ is full by Theorem \ref{definab}(ii) and faithful by the conservativity of $M_{\mathbb T}$. Thus, being essentially surjective on ${\textrm{f.p.} {\mathbb T}\textrm{-mod}(\Set)}_{s}$, $d$ gives an equivalence of categories $\tilde{{\cal C}_{\mathbb T}}\simeq {\textrm{f.p.} {\mathbb T}\textrm{-mod}(\Set)}_{s}^{\textrm{op}}$, as required.          
\end{proofs}

The theorem motivates the following definition: given a theory of presheaf type $\mathbb T$, we say that the category $\textrm{f.p.} {\mathbb T}\textrm{-mod}(\Set)$ is syntactically presented if every finitely presentable model in $\Set$ is finitely presented; in fact, $\textrm{f.p.} {\mathbb T}\textrm{-mod}(\Set)$ is syntactically presented if and only if the functor $d$ yields an equivalence of categories between $\tilde{{\cal C}_{\mathbb T}}$ and $\textrm{f.p.}{\mathbb T}\textrm{-mod}(\Set)^{\textrm{op}}$.\\    
For fragments of geometric logic such as cartesian, regular or coherent logic, it also makes sense to define the notion of finitely presented model; indeed, one can say that a given model in $\Set$ of a cartesian (resp. regular, coherent) theory $\mathbb T$ is finitely presented if there exists a cartesian (resp. regular, coherent) formula over the signature of $\mathbb T$ which presents it. In fact, from \cite{El2} we know that all the finitely presented models of a cartesian theory (in the geometric sense) are presented by cartesian formulae. From the proof of the theorem and Remark \ref{restr}, it is clear that the analogue of Theorem \ref{dualitap} for cartesian, regular or coherent logic also holds.\\
The following result shows that arrows in syntactic categories behave naturally with respect to different fragments of logic i.e. for any cartesian (resp. regular, coherent) theory $\mathbb T$, the cartesian (resp. regular, coherent) syntactic category of $\mathbb T$ is a full subcategory of the geometric syntactic category of $\mathbb T$.

\begin{theorem}
Let $\mathbb T$ be a cartesian (resp. regular, coherent) theory. If $[\theta]:\{\vec{y}. \psi\}\to \{\vec{x}. \phi\}$ is an arrow in the geometric syntactic category ${\cal C}_{\mathbb T}$ of $\mathbb T$ and both $\{\vec{y}. \psi\}$ and $\{\vec{x}. \phi\}$ are cartesian (resp. regular, coherent) formulae then $\theta(\vec{y}, \vec{x})$ is $\mathbb T$-provably equivalent to a cartesian (resp. regular, coherent) formula $\theta'(\vec{y}, \vec{x})$ which is $\mathbb T$-provably functional from $\{\vec{y}. \psi\}$ to $\{\vec{x}. \phi\}$ i.e. which gives an arrow $\{\vec{y}. \psi\}\to \{\vec{x}. \phi\}$ in the cartesian (resp. regular, coherent) syntactic category ${\cal C}_{\mathbb T}^{\textrm{cart}}$ of $\mathbb T$.
\end{theorem}

\begin{proofs}
We can deduce the result from the logical completeness of universal models as follows. Given the Yoneda embedding $y:{\cal C}_{\mathbb T}\to \Sh({\cal C}_{\mathbb T}, J_{\mathbb T})$, by the cartesian (resp. regular, coherent) version of Theorem \ref{definab}(ii), $y([\theta])=[[\theta]]_{U_{\mathbb T}}$ is equal to the interpretation in the universal model $U_{\mathbb T}$ of $\mathbb T$ in $\Sh({\cal C}_{\mathbb T}, J_{\mathbb T}) \simeq [({\cal C}^{\textrm{cart}}_{\mathbb T})^{\textrm{op}}, \Set]$ (resp. $\Sh({\cal C}_{\mathbb T}, J_{\mathbb T}) \simeq \Sh({\cal C}^{\textrm{reg}}_{\mathbb T}, J_{\mathbb T}^{\textrm{reg}})$, $\Sh({\cal C}_{\mathbb T}, J_{\mathbb T}) \simeq \Sh({\cal C}^{\textrm{coh}}_{\mathbb T}, J_{\mathbb T}^{\textrm{coh}})$) of a cartesian (resp. regular, coherent) formula $\theta'$; then $\theta$ and $\theta'$ have the same interpretation in the universal model $U_{\mathbb T}$ and hence, by the conservativity of universal models, they are $\mathbb T$-provably equivalent.   \end{proofs}

\section{Universal models of quotients of a theory of presheaf type}\label{terza}
Suppose we are given a theory of presheaf type $\mathbb T$ over a signature $\Sigma$ and a quotient ${\mathbb T}'$ of $\mathbb T$ obtained from $\mathbb T$ by adding axioms $\sigma$ of the form $\phi \vdash_{\vec{x}} \mathbin{\mathop{\textrm{\huge $\vee$}}\limits_{i\in I}}(\exists \vec{y_{i}})\theta_{i}$, where, for any $i\in I$, $[\theta_{i}]:\{\vec{y_{i}}. \psi\}\to \{\vec{x}. \phi\}$ is an arrow in ${\cal C}_{\mathbb T}$ and $\phi(\vec{x})$, $\psi(\vec{y_{i}})$ are geometric formulae over $\Sigma$ presenting respectively ${\mathbb T}$-models $M_{\phi}$ and $M_{\psi_{i}}$. Denote by $S_{\sigma}$ the cosieve in $\textrm{f.p.} {\mathbb T}\textrm{-mod}(\Set)$ on $M_{\phi}$ generated by the arrows $d([\theta_{i}])$ as $i$ varies in $I$. Then we have the following result.
\begin{theorem}\label{constru}
Let $\mathbb T$ be a theory of presheaf type over a signature $\Sigma$ and ${\mathbb T}'$ a quotient of $\mathbb T$ obtained from $\mathbb T$ by adding axioms $\sigma$ of the form $\phi \vdash_{\vec{x}} \mathbin{\mathop{\textrm{\huge $\vee$}}\limits_{i\in I}}(\exists \vec{y_{i}})\theta_{i}$, where, for each $i\in I$, $[\theta_{i}]:\{\vec{y_{i}}. \psi\}\to \{\vec{x}. \phi\}$ is an arrow in ${\cal C}_{\mathbb T}$ and $\phi(\vec{x})$, $\psi(\vec{y_{i}})$ are geometric formulae over $\Sigma$ presenting respectively ${\mathbb T}$-models $M_{\phi}$ and $M_{\psi_{i}}$. With the notation above, ${\mathbb T}'$ is classified by the topos $\Sh(\textrm{f.p.} {\mathbb T}\textrm{-mod}(\Set)^{\textrm{op}}, J)$ where $J$ is the Grothendieck topology on $\textrm{f.p.} {\mathbb T}\textrm{-mod}(\Set)^{\textrm{op}}$ generated by the sieves $S_{\sigma}$. 
\end{theorem} 

\begin{proofs}
We know from \cite{OC6} that ${\mathbb T}'$ corresponds to a unique Grothendieck topology $J$ on $\textrm{f.p.} {\mathbb T}\textrm{-mod}(\Set)^{\textrm{op}}$ such that the given equivalence ${\mathbb T}\textrm{-mod}({\cal E}) \simeq {\bf Flat}(\textrm{f.p.} {\mathbb T}\textrm{-mod}(\Set)^{\textrm{op}}, {\cal E})$ restricts to an equivalence ${\mathbb T}'\textrm{-mod}({\cal E}) \simeq {\bf Flat}_{J}(\textrm{f.p.} {\mathbb T}\textrm{-mod}(\Set)^{\textrm{op}}, {\cal E})$ (naturally in ${\cal E}\in \mathfrak{BTop}$). Indeed, this equivalence is induced by transferring the subtopos of $\Sh({\cal C}_{J_{\mathbb T}}, J_{{\mathbb T}'}^{\mathbb T})$ of $\Sh({\cal C}_{\mathbb T}, J_{\mathbb T})$ correponding to ${\mathbb T}'$ via Theorem 3.6 \cite{OC6} to a geometric inclusion $\Sh(\textrm{f.p.} {\mathbb T}\textrm{-mod}(\Set)^{\textrm{op}}, J)\hookrightarrow [\textrm{f.p.} {\mathbb T}\textrm{-mod}(\Set), \Set]$ via the equivalence $[\textrm{f.p.} {\mathbb T}\textrm{-mod}(\Set), \Set]\simeq \Sh({\cal C}_{\mathbb T}, J_{\mathbb T})$ of classifying toposes of $\mathbb T$. We want to prove that $J$ is the Grothendieck topology generated by the sieves $S_{\sigma}$. From section 14.2 \cite{OC6} we know that for each axiom $\sigma$ and functor $F\in {\bf Flat}(\textrm{f.p.} {\mathbb T}\textrm{-mod}(\Set)^{\textrm{op}}, {\cal E})$, $F$ sends $S_{\sigma}$ to an epimorphic family if and only if $\sigma$ holds in the corresponding $\mathbb T$-model. So from the equivalence ${\mathbb T}'\textrm{-mod}({\cal E}) \simeq {\bf Flat}_{J}(\textrm{f.p.} {\mathbb T}\textrm{-mod}(\Set)^{\textrm{op}}, {\cal E})$ we deduce that for any ${\cal E}\in \mathfrak{BTop}$ and any $F\in {\bf Flat}_{J}(\textrm{f.p.} {\mathbb T}\textrm{-mod}(\Set)^{\textrm{op}}, {\cal E})$, $F$ sends $S_{\sigma}$ to an epimorphic family. This implies, by Lemma 3 p. 393 \cite{MM} applied to the canonical geometric inclusion $\Sh(\textrm{f.p.} {\mathbb T}\textrm{-mod}(\Set)^{\textrm{op}}, J)\hookrightarrow [\textrm{f.p.} {\mathbb T}\textrm{-mod}(\Set), \Set]$, that the associated sheaf functor $a_{J}:[\textrm{f.p.} {\mathbb T}\textrm{-mod}(\Set), \Set] \to \Sh(\textrm{f.p.} {\mathbb T}\textrm{-mod}(\Set)^{\textrm{op}}, J)$ sends the monomorphism $S_{\sigma}\mono M_{\phi}$ to an isomorphism i.e. $S_{\sigma}$ is $J$-covering. Alternatively, this follows from Theorem \ref{dualint} and Lemma \ref{subtopos} since $a_{J}$ preserves the interpretation of geometric formulae. Thus the Grothendieck topology $J'$ generated by the $S_{\sigma}$ is contained in $J$. To prove that $J'=J$, it is equivalent to verify, by Diaconescu's theorem and the 2-dimensional Yoneda Lemma, that for any Grothendieck topos $\cal E$, ${\bf Flat}_{J'}(\textrm{f.p.} {\mathbb T}\textrm{-mod}(\Set)^{\textrm{op}}, {\cal E})={\bf Flat}_{J}(\textrm{f.p.} {\mathbb T}\textrm{-mod}(\Set)^{\textrm{op}}, {\cal E})$. Since $J'\subseteq J$, clearly ${\bf Flat}_{J}(\textrm{f.p.} {\mathbb T}\textrm{-mod}(\Set)^{\textrm{op}}, {\cal E}) \subseteq {\bf Flat}_{J'}(\textrm{f.p.} {\mathbb T}\textrm{-mod}(\Set)^{\textrm{op}}, {\cal E})$. To prove the other inclusion, observe that, since $J'$ contains all the sieves $S_{\sigma}$ then the $\mathbb T$-model corresponding to a functor $F\in {\bf Flat}_{J'}(\textrm{f.p.} {\mathbb T}\textrm{-mod}(\Set)^{\textrm{op}}, {\cal E})$ belongs to ${\mathbb T}'\textrm{-mod}({\cal E})$ and hence, by the equivalence ${\mathbb T}'\textrm{-mod}({\cal E}) \simeq {\bf Flat}_{J}(\textrm{f.p.} {\mathbb T}\textrm{-mod}(\Set)^{\textrm{op}}, {\cal E})$, $F$ is $J$-continuous, as required.    
\end{proofs}

Below, we will refer to the Grothendieck topology on $\textrm{f.p.} {\mathbb T}\textrm{-mod}(\Set)^{\textrm{op}}$ corresponding to a quotient ${\mathbb T}'$ of $\mathbb T$ as in the proof of Theorem \ref{constru} as the associated $\mathbb T$-topology of ${\mathbb T}'$.\\

In passing, we observe the following fact.

\begin{proposition}
Let $\mathbb T$ be a theory of presheaf type over a signature $\Sigma$ and ${\mathbb T}'$ a quotient of $\mathbb T$ with associated $\mathbb T$-topology $J$ on $\textrm{f.p.} {\mathbb T}\textrm{-mod}(\Set)^{\textrm{op}}$. Then, for any geometric sequent $\sigma$ over $\Sigma$ of the form $\phi \vdash_{\vec{x}} \mathbin{\mathop{\textrm{\huge $\vee$}}\limits_{i\in I}}(\exists \vec{y_{i}})\theta_{i}$, where, for each $i\in I$, $[\theta_{i}]:\{\vec{y_{i}}. \psi\}\to \{\vec{x}. \phi\}$ is an arrow in ${\cal C}_{\mathbb T}$ and $\phi(\vec{x})$, $\psi(\vec{y_{i}})$ are geometric formulae over $\Sigma$ presenting respectively ${\mathbb T}$-models $M_{\phi}$ and $M_{\psi_{i}}$, $\sigma$ is provable in ${\mathbb T}'$ if and only if the sieve $S_{\sigma}$ is $J$-covering.    
\end{proposition}

\begin{proofs}
This immediately follows from Theorem \ref{dualint} and Lemma \ref{subtopos} by using the fact that $a_{J}$ preserves the interpretation of geometric formulae and the conservativity of universal models. 
\end{proofs} 

We have seen that if $J$ is the associated $\mathbb T$-topology of a quotient ${\mathbb T}'$ of $\mathbb T$ then (by Lemma \ref{subtopos}) the $\Sigma$-structure $M_{\mathbb T}^{J}:=a_{J}(M_{\mathbb T})$ is a universal model of ${\mathbb T}'$, where $a_{J}:[\textrm{f.p.} {\mathbb T}\textrm{-mod}(\Set), \Set]\to \Sh(\textrm{f.p.} {\mathbb T}\textrm{-mod}(\Set), J)$ is the associated sheaf functor and $M_{\mathbb T}$ is the universal model of $\mathbb T$ as in\\ Theorem \ref{presheaf}.\\
Let us now suppose that that for any sort $A$ in the signature of $\mathbb T$ the formula $\{x^{A}.\top\}$ presents a $\mathbb T$-model and the topology $J$ is subcanonical - we will refer to the conjunction of these conditions as $(\ast)$ - and derive, under these assumptions, a simple description of the universal model of ${\mathbb T}'$.\\
First, we note that $(a_{J}(M_{\mathbb T}))A=M_{\mathbb T}A$ for any sort $A$ and $(a_{J}(M_{\mathbb T}))f=M_{\mathbb T}f$ for any function symbol $f$. Next, we observe that if $\cal E$ is an elementary topos and $j$ a local operator on $\cal E$ corresponding to a universal closure operator $c_{j}$ on $\cal E$ then, given a subobject $m:A'\mono A$ in $\cal E$, $a_{j}(m)\cong a_{J}(c_{j}(m))$; indeed, $c_{j}(m)$ is defined by the following pullback square
\[  
\xymatrix {
c_{j}(A') \ar[r] \ar[d] & a_{j}A' \ar[d]^{a_{j}m}\\
A \ar[r]^{\eta_{A}^{a_{j}}} & a_{j}A }
\] 
and $a_{j}$ preserves pullbacks. In particular, if $A$ is a $j$-sheaf then $a_{j}(m)=c_{j}(m)$.\\
Thus, under $(\ast)$, given a geometric formula $\phi(\vec{x})$ over $\Sigma$, where $\vec{x}=(x^{A_{1}}, \ldots, x^{A_{n}})$, the interpretation $[[\vec{x}. \phi]]_{M_{\mathbb T}^{J}}$ of $\phi$ in the universal model $M_{\mathbb T}^{J}$ is equal to the $J$-closure of $[[\vec{x}. \phi]]_{M_{\mathbb T}}\mono M_{\mathbb T}A_{1}\times \cdots \times M_{\mathbb T}A_{n}$ in $[\textrm{f.p.} {\mathbb T}\textrm{-mod}(\Set), \Set]$.\\
Given $\phi=\phi(x^{A_{1}}, \cdots, x^{A_{n}})$ and $\vec{a}=(a_{1}, \cdots, a_{n})\in MA_{1}\times \cdots MA_{n}$ where $M\in \textrm{f.p.} {\mathbb T}\textrm{-mod}(\Set)$, let us define the sieve\\
$S^{\phi}_{\vec{a}}:=\{f:M\to N \textrm{ in f.p.} {\mathbb T}\textrm{-mod}(\Set) \textrm{ | } f(\vec{a})\in [[\vec{x}.\phi]]_{N}\}$\\
in ${\textrm{f.p.} {\mathbb T}\textrm{-mod}(\Set)}^{\textrm{op}}$.\\
Then, by applying formula (6) p. 235 \cite{MM}, we obtain that the $J$-closure $\overline{[[\vec{x}. \phi]]_{M_{\mathbb T}}}^{J}$ of $[[\vec{x}. \phi]]_{M_{\mathbb T}}\mono M_{\mathbb T}A_{1}\times \cdots \times M_{\mathbb T}A_{n}$ in $[\textrm{f.p.} {\mathbb T}\textrm{-mod}(\Set), \Set]$ is given by:       
\[
\begin{array}{ccl}
\overline{[[\vec{x}. \phi]]_{M_{\mathbb T}}}^{J}(M) & = & \{(a_{1}, \cdots, a_{n})\in MA_{1}\times \cdots MA_{n} \textrm{ | } S^{\phi}_{\vec{a}}\in J(M)\}   
\end{array}
\]
for any $M\in \textrm{f.p.} {\mathbb T}\textrm{-mod}(\Set)$. 
Accordingly, we say that $\phi$ is $J$-closed if $\overline{[[\vec{x}. \phi]]_{M_{\mathbb T}}}^{J}=[[\vec{x}. \phi]]_{M_{\mathbb T}}$ i.e. for any $\vec{a}\in MA_{1}\times \cdots MA_{n}$, $\{f:M\to N \textrm{ in f.p.} {\mathbb T}\textrm{-mod}(\Set) \textrm{ | } f(\vec{a})\in [[\phi]]_{N}\}\in J(M)$ implies $\vec{a}\in [[\phi]]_{M}$.\\
\begin{eg}
\emph{Let us consider the theory $\mathbb T$ of commutative rings with unit. Then, by the proof of Lemma 6.3.(ii) \cite{OC3}, the formula $\phi(x)=(\exists y)(x\cdot y=1)$ is $J$-closed where $J$ is the Zariski topology on ${\textrm{f.p.} {\mathbb T}\textrm{-mod}(\Set)}^{\textrm{op}}$ (cfr. also section \ref{exa} below).}   
\end{eg}

Hence, the conservativity of universal models yields the following result.

\begin{theorem}\label{subcan}
Let ${\mathbb T}'$ be a quotient of a theory of presheaf type $\mathbb T$ over a signature $\Sigma$, with associated $\mathbb T$-topology $J$ on $\textrm{f.p.} {\mathbb T}\textrm{-mod}(\Set)^{\textrm{op}}$. Given a geometric sequent $\phi \vdash_{\vec{x}} \psi$ over $\Sigma$, if for any $M\in \textrm{f.p.} {\mathbb T}\textrm{-mod}(\Set)$ and any $\vec{a}\in M$, $S^{\phi}_{\vec{a}}\in J(M)$ implies $S^{\psi}_{\vec{a}}\in J(M)$ then $\phi \vdash_{\vec{x}} \psi$ is provable in ${\mathbb T}'$. The converse holds if for each sort $A$ of $\Sigma$, $\{x^{A}. \top\}$ presents a $\mathbb T$-model and $J$ is subcanonical.\qed   
\end{theorem}

We recall from that a geometric theory $\mathbb{T}$ over a signature $\Sigma$ is said to have enough models if for every geometric sequent $\sigma$ over $\Sigma$, $M\vDash \sigma$ for all the $\mathbb T$-models $M$ in $\Set$ implies that $\sigma$ is provable in $\mathbb T$ using geometric logic.\\  

The following result is an application of universal models to determine whether a quotient of a theory of presheaf type $\mathbb T$ has enough models.

\begin{corollary}
Let ${\mathbb T}'$ be a quotient of a theory of presheaf type $\mathbb T$ with associated $\mathbb T$-topology $J$ on $\textrm{f.p.} {\mathbb T}\textrm{-mod}(\Set)^{\textrm{op}}$. If for any $M\in \textrm{f.p.} {\mathbb T}\textrm{-mod}(\Set)$ there exists a $J$-covering cosieve in $\textrm{f.p.} {\mathbb T}\textrm{-mod}(\Set)$ on $M$ generated by arrows in $\textrm{f.p.} {\mathbb T}\textrm{-mod}(\Set)$ whose codomain is a ${\mathbb T}'$-model then ${\mathbb T}'$ has enough models.   
\end{corollary}
\begin{proofs}
In view of Theorem \ref{subcan}, it is enough to prove that if for any $M\in \textrm{f.p.} {\mathbb T}\textrm{-mod}(\Set)$ there exists a $J$-covering cosieve $S_{M}$ in $\textrm{f.p.} {\mathbb T}\textrm{-mod}(\Set)^{\textrm{op}}$ on $M$ generated by arrows in $\textrm{f.p.} {\mathbb T}\textrm{-mod}(\Set)$ whose codomain is a ${\mathbb T}'$-model then for any geometric sequent $\phi \vdash_{\vec{x}} \psi$ over $\Sigma$ with $\vec{x}=(x^{A_{1}}, \ldots, x^{A_{n}})$ which is satisfied in any ${\mathbb T}'$-model in $\Set$, for each $M\in \textrm{f.p.} {\mathbb T}\textrm{-mod}(\Set)$ and $\vec{a}\in MA_{1}\times \cdots \times MA_{n}$, $S^{\phi}_{\vec{a}}\in J(M)$ implies $S^{\psi}_{\vec{a}}\in J(M)$. Now, if $S_{M}\in J(M)$ then $S^{\phi}_{\vec{a}}\in J(M)$ implies that $S^{\phi}_{\vec{a}}\cap S_{M}\in J(M)$; but $S^{\phi}_{\vec{a}}\cap S_{M}\subseteq S^{\psi}_{\vec{a}}\cap S_{M}\subseteq S^{\psi}_{\vec{a}}$, and hence $S^{\psi}_{\vec{a}}\in J(M)$, as required.      
\end{proofs} 

In view of the results above, it is useful to find criteria for a Grothendieck topology $J$ on a category $\cal C$ to be subcanonical.\\
In fact, from the discussion after the proof of Corollary C2.1.11 \cite{El2}, we can extract the following criterion: $J$ is subcanonical if and only if all its $J$-covering sieves are effective-epimorphic. If $\cal F$ is a collection of sieves in $\cal C$ which generates $J$ then it follows from Proposition C2.1.9 \cite{El2} that $J$ is subcanonical if and only if every sieve in $\cal F$ is universally effective-epimorphic; in particular, if $\cal F$ is closed under sieve pullbacks then $J$ is subcanonical if and only if every sieve in $\cal F$ is effective-epimorphic.\\
In the case $\textrm{f.p.} {\mathbb T}\textrm{-mod}(\Set)$ is syntactically presented, we can rephrase in terms of definable subsets the condition for a sieve $\{f_{i}:M\to M_{i} \textrm{ in $\textrm{f.p.} {\mathbb T}\textrm{-mod}(\Set)$ | } i\in I\}$ in $\textrm{f.p.} {\mathbb T}\textrm{-mod}(\Set)^{\textrm{op}}$ to be effective-epimorphic, as follows: for any formula $\{\vec{z}. \chi\}$ in $\tilde{{\cal C}_{\mathbb T}}$ and any collection $\{c_{i}\in [[\vec{z}.\chi]]_{M_{i}} \textrm{ | } i\in I\}$ of elements such that for any arrow $f:M_{i}\to M_{j}$ in $\textrm{f.p.} {\mathbb T}\textrm{-mod}(\Set)$ with the property that $f\circ f_{i}=f_{j}$, $f(c_{i})=c_{j}$, there exists a unique $c\in [[\vec{z}. \chi]]_{M}$ such that for any $i\in I$, $f_{i}(c)=c_{i}$.

\section{Applications to definability}

In this section, we exhibit a link between notions of definability of properties of elements of finitely presentable models of a given theory of presheaf type $\mathbb T$ and properties of the universal models of quotients of $\mathbb T$.\\
Let us start with some general remarks.\\

Let $\mathbb T$ be a geometric theory over a signature $\Sigma$ and $M$ a $\mathbb T$-model in a Grothendieck topos $\cal E$. For each context $\vec{x}=(x_{1}^{A_{1}}, \ldots, x_{n}^{A_{n}})$ over $\Sigma$, the subobjects of $MA_{1}\times \ldots \times MA_{n}$ of the form $[[\vec{x}. \phi]]_{M}$, where $\phi(\vec{x})$ is a geometric formula over $\Sigma$ in the context $\vec{x}$, clearly form a subframe of $\Sub_{\cal E}(MA_{1}\times \ldots \times MA_{n})$, which we denote by $\textrm{Def}_{\vec{x}}^{\textrm{geom}}(M)$. Note that there is a geometric surjective functor $\textrm{Int}_{\vec{x}}^{M}:\Sub_{{\cal C}_{\mathbb T}}(\{\vec{x}.\top\})\rightarrow \textrm{Def}_{\vec{x}}^{\textrm{geom}}(M)$ which sends each formula $\phi(\vec{x})$ (identified with the corresponding subobject $[\phi]:\{\vec{x}.\phi\}\rightarrow \{\vec{x}.\top\}$ in ${\cal C}_{\mathbb T}$ via Lemma D1.4.4(iv) \cite{El2}) to its interpretation $[[\vec{x}. \phi]]_{M}$ in the model $M$, and if $M$ is conservative then this functor is an isomorphism; in particular, the pseudocomplementation in $\textrm{Def}_{\vec{x}}^{\textrm{geom}}(M)$ of an object $[[\vec{x}. \phi]]_{M} \in \textrm{Def}_{\vec{x}}^{\textrm{geom}}(M)$ coincides with $[[\vec{x}. \neg^{\mathbb T} \phi]]_{M}$, where $[\neg^{\mathbb T}\phi]$ is the pseudocomplementation of $[\phi]$ in $\Sub_{{\cal C}_{\mathbb T}}(\{\vec{x}. \top\})$. If moreover any subobject $S\mono MA_{1}\times \cdots MA_{n}$ in $\cal E$ is the interpretation in $M$ of a geometric formula $\phi(x^{A_{1}}, \ldots, x^{A_{n}})$ over $\Sigma$ then $\textrm{Def}_{\vec{x}}^{\textrm{ geom}}(M)=\Sub_{\cal E}(MA_{1}\times \cdots MA_{n})$ and hence the pseudocomplementation of an object $[[\vec{x}. \phi]]_{M}$ in $\textrm{Def}_{\vec{x}}^{\textrm{geom}}(M)$ coincides with its pseudocomplementation $\neg [[\vec{x}. \phi]]_{M}$ in the subobject lattice $\Sub_{\cal E}(MA_{1}\times \cdots MA_{n})$. Thus, since universal models are always conservative and satisfy the condition of Theorem \ref{definab}(i), we get the following result.

\begin{theorem}\label{unif}
Let $\mathbb T$ be a geometric theory and $M$ a universal model of $\mathbb T$ in a Grothendieck topos $\cal E$. Then for any context $\vec{x}=(x^{A_{1}}, \ldots, x^{A_{n}})$, we have Heyting algebra isomorphisms
\[
\textrm{Int}_{\vec{x}}^{M}: \Sub_{{\cal C}_{\mathbb T}}(\{\vec{x}.\top\}) \stackrel{\cong}{\longrightarrow} \textrm{Def}_{\vec{x}}^{\textrm{ geom}}(M)=\Sub_{\cal E}(MA_{1}\times \cdots MA_{n})
\]
\end{theorem}\qed

The following lemma concerns the interpretation of pseudocomplements $\neg^{\mathbb T} \phi$ in models of $\mathbb T$; below, the symbol $\neg$ denotes the operation of pseudocomplementation of subobjects in the topos in which the relevant $\mathbb T$-model lies.\\
 
\begin{lemma}\label{lemmafond}
Let $\mathbb T$ be a geometric theory over a signature $\Sigma$ and $M$ a $\mathbb T$-model in a topos $\cal E$. With the above notation, we have:\\
(i) For any geometric formula $\phi$ over $\Sigma$, $[[\neg^{\mathbb T}\phi]]_{M}\leq \neg [[\phi]]_{M}$;\\   
(ii) For any geometric formulae $\phi_{1}$ and $\phi_{2}$ over $\Sigma$ in the same context such that $\phi_{1}\wedge \phi_{2}\vdash \bot$ is provable in $\mathbb T$, if $\neg[[\phi_{1}]]_{M}=[[\phi_{2}]]_{M}$ then $[[\neg^{\mathbb T} \phi_{1}]]_{M}=[[\phi_{2}]]_{M}$. In particular, if $M$ is conservative then $\neg^{\mathbb T} \phi_{1}$ and $\phi_{2}$ are $\mathbb T$-provably equivalent;\\
(iii) If $M$ is conservative then, given a geometric formula $\phi$ over $\Sigma$, $[[\neg^{\mathbb T} \phi]]_{M}=\neg [[\phi]]_{M}$ if and only if $\neg [[\phi]]_{M}$ is the interpretation in $M$ of a geometric formula over $\Sigma$. In particular, if $M$ is a universal model of $\mathbb T$ then $[[\neg^{\mathbb T} \phi]]_{M}=\neg [[\phi]]_{M}$ for any geometric formula $\phi$ over $\Sigma$.  
\end{lemma}
\begin{proofs}
(i) Since $\neg^{\mathbb T}\phi \wedge \phi$ is $\mathbb T$-provably equivalent to $\top$ then $[[\neg^{\mathbb T}\phi]]_{M}\cap[[\phi]]_{M}=0$ i.e. $[[\neg^{\mathbb T}\phi]]_{M}\leq \neg [[\phi]]_{M}$.\\ 
(ii) Since $\neg[[\phi_{1}]]_{M}=[[\phi_{2}]]_{M}$, by part (i) $[[\neg^{\mathbb T} \phi_{1}]]_{M}\leq \neg[[\phi_{1}]]_{M}$ and hence $\neg[[\phi_{1}]]_{M}=[[\phi_{2}]]_{M}$ implies $[[\neg^{\mathbb T} \phi_{1}]]_{M}\leq [[\phi_{2}]]_{M}$; the converse inequality follows from the fact that $\{\vec{x}.\phi_{2}\}\leq \{\vec{x}. \neg \phi_{1}\}$ in $\Sub_{{\cal C}_{\mathbb T}}(\{\vec{x}.\top\})$.\\
(iii) The `only if' direction is obvious. To prove the other direction, we observe that if $\neg [[\phi]]_{M}=[[\chi]]_{M}$ then $[[\phi]]_{M}\cap[[\chi]]_{M}=0$, which implies $\phi \wedge \chi \vdash \bot$ provable in $\mathbb T$, since $M$ is conservative; thus the thesis follows from part (ii) of the lemma. The last assertion now follows from Theorem \ref{definab}(i).    
\end{proofs}

The lemma can be profitably applied in the context of theories of presheaf type, where we have an explicit description of their universal models (Theorem \ref{presheaf}). To this end, let us recall the explicit description of the pseudocomplementation of subobjects in presheaf toposes given from formula (16) p. 273 \cite{MM}.\\
The pseudocomplementation $\neg A \mono E$ of a subobject $A\mono E$ in a presheaf topos $[{\cal C}^{\textrm{op}}, \Set]$ is given by:\\ 
\[
(\neg A)(c):=\{e\in E(c) \textrm{ for all } f:d\to c \textrm{ in $\cal C$, } E(f)(e)\notin A(d)\}.
\]
Note in particular that, given a sieve $R$ on an object $c\in {\cal C}$, regarded as a subobject $R\mono {\cal C}(-, c)$ in $[{\cal C}^{\textrm{op}}, \Set]$, we have\\
\[
\neg R:=\{f:d\to c \textrm{ in } {\cal C} \textrm{ | } f^{\ast}(R)=\emptyset\},
\]
and
\[
\neg \neg R:=\{f:d\to c \textrm{ in } {\cal C} \textrm{ | } f^{\ast}(R) \textrm{ is stably non-empty}\}.
\]

\begin{proposition}\label{pseudoneg}
Let $\mathbb T$ be a theory of presheaf type over a signature $\Sigma$ and $M$ a finitely presentable $\mathbb T$-model in $\Set$. Given a geometric formula $\phi(x^{A_{1}}, \ldots, x^{A_{n}})$ over $\Sigma$,
\[
\begin{array}{ccl}
[[\neg^{\mathbb T} \phi]]_{M} & = & \{(a_{1}, \cdots, a_{n})\in MA_{1}\times \cdots MA_{n} \textrm{ | for all } f:M\to N \\
& & \textrm{ in f.p.} {\mathbb T}\textrm{-mod}(\Set) \textrm{ | } f(a_{1}, \cdots a_{n})\notin [[\phi]]_{N}\}   
\end{array}
\]
\end{proposition}

\begin{proofs}
By Lemma \ref{lemmafond}, $[[\neg \phi]]_{M_{\mathbb T}}=\neg [[\phi]]_{M_{\mathbb T}}$, where $M_{\mathbb T}$ is the universal model of $\mathbb T$ in $[\textrm{f.p.} {\mathbb T}\textrm{-mod}(\Set), \Set]$; then the thesis follows from the above-mentioned formula for pseudocomplementation in presheaf toposes and Theorem \ref{presheaf}.
\end{proofs}

Let $\mathbb T$ be a theory of presheaf type over a signature $\Sigma$ and $\phi(\vec{x})$ a geometric formula over $\Sigma$, where $\vec{x}=(x^{A_{1}}, \ldots, x^{A_{n}})$.\\
One can consider the property $P_{\{\vec{x}. \phi\}}$ (of elements of finitely presentable $\mathbb T$-models) defined by:\\
for any $M\in \textrm{f.p.} {\mathbb T}\textrm{-mod}(\Set)$ and $\vec{a}\in MA_{1}\times \cdots \times MA_{n}$, $M \vDash P_{\phi}(\vec{a})$ if and only if $M \nvDash \phi(\vec{a})$,\\
and ask when $P_{\{\vec{x}. \phi\}}$ is definable by a geometric formula over $\Sigma$.\\ 
From Lemma \ref{lemmafond} it follows that for any $\phi(\vec{x})$, $P_{\{\vec{x}. \phi\}}$ is definable if and only if it is definable by $\neg^{\mathbb T} \phi$, if and only if for any $M\in \textrm{f.p.} {\mathbb T}\textrm{-mod}(\Set)$ and $\vec{a}\in MA_{1}\times \cdots \times MA_{n}$, $M \nvDash \phi(\vec{a})$ if and only if for all $f:M\to N$ in $\textrm{f.p.} {\mathbb T}\textrm{-mod}(\Set)$, $N \nvDash \phi(\vec{a})$ (note that this characterization also follows from Corollary \ref{predef}).\\
Also, one can consider the property $Q^{{\mathbb T}'}_{\{\vec{x}.\phi\}}:=P_{\{\vec{x}.\neg^{{\mathbb T}'} \phi\}}$, where ${\mathbb T}'$ is any quotient of $\mathbb T$.\\
From Proposition \ref{pseudoneg}, we get the following explicit description of $Q^{\mathbb T}_{\{\vec{x}.\phi\}}$: for any $M\in \textrm{f.p.} {\mathbb T}\textrm{-mod}(\Set)$ and $\vec{a}\in MA_{1}\times \cdots \times MA_{n}$, $M \vDash Q^{\mathbb T}_{\{\vec{x}.\phi\}}(\vec{a})$ if and only if there exists $f:M\to N$ in $\textrm{f.p.} {\mathbb T}\textrm{-mod}(\Set)$ such that $f(\vec{a})\in [[\phi]]_{N}$. As we shall see below, this property plays an important role in connection with De Morgan's law on Grothendieck toposes.\\
From Lemma \ref{lemmafond} it follows that for any $\phi(\vec{x})$, $Q^{\mathbb T}_{\{\vec{x}.\phi\}}$ is definable by a geometric formula over $\Sigma$ if and only if it is definable by $\neg^{\mathbb T} \neg^{\mathbb T} \phi$, if and only for any $M\in \textrm{f.p.} {\mathbb T}\textrm{-mod}(\Set)$ and $\vec{a}\in MA_{1}\times \cdots \times MA_{n}$, $M \nvDash Q^{\mathbb T}_{\{\vec{x}.\phi\}}(\vec{a})$ if and only if $\vec{a}\in \neg\neg [[\phi]]_{M}$, if and only if (by formula (15) p. 273 \cite{MM}) for all $f:M\to N$ in $\textrm{f.p.} {\mathbb T}\textrm{-mod}(\Set)$ there is $g:N\to P$ in in $\textrm{f.p.} {\mathbb T}\textrm{-mod}(\Set)$ such that $(g\circ f)(\vec{a})\in [[\phi]]_{P}$.\\ 
Another natural property to consider is $Z^{{\mathbb T}'}_{\{\vec{x}.\phi\}}:=P_{\{\vec{x}. \phi\vee \neg^{{\mathbb T}'} \phi\}}$, where ${\mathbb T}'$ is a quotient of $\mathbb T$. In connection with this, we note that in order to establish that a certain formula $\chi$ is of the form $\neg^{{\mathbb T}'} \phi$, it is convenient to work in a universal model $U$ of ${\mathbb T}'$ (as it is done for example in \cite{OC4}) and prove that its interpretation here coincides with $\neg [[\phi]]_{U}$ (cfr. Lemma \ref{lemmafond}).\\  
       
We can now fulfil a promise made at the beginning of this section.\\
Below, the relation of $\mathbb T$-provable equivalence of geometric formulae in the same context is denoted by $\stackrel{\mathbb T}{\sim}$. 

\begin{theorem}\label{teofond}
Let $\mathbb T$ be a theory of presheaf type $\mathbb T$ over a signature $\Sigma$, ${\mathbb T}'$ a quotient of $\mathbb T$, $\phi(\vec{x})$ a geometric formula over $\Sigma$, and let $P$ be either property $P_{\{\vec{x}. \phi\}}$ or property $Z^{{\mathbb T}'}_{\{\vec{x}.\phi\}}$. Then\\
(i) If $P$ is definable by a geometric formula over $\Sigma$ then $\phi \vee \neg^{{\mathbb T}'} \phi \stackrel{{\mathbb T}'}{\sim} \top$;\\
(ii) The converse implication holds if ${\mathbb T}'={\mathbb T}$.
\end{theorem}
\begin{proofs}
(i) By Lemma \ref{subtopos}, $U':=a_{J}(M_{\mathbb T})$ is a universal model of ${\mathbb T}'$ in ${\cal F}:=\Sh(\textrm{f.p.} {\mathbb T}\textrm{-mod}(\Set)^{\textrm{op}}, J)$, where $J$ is the associated $\mathbb T$-topology of ${\mathbb T}'$, $M_{\mathbb T}$ is the universal model of $\mathbb T$ in ${\cal E}:=[\textrm{f.p.} {\mathbb T}\textrm{-mod}(\Set), \Set]$ and $a_{J}:[\textrm{f.p.} {\mathbb T}\textrm{-mod}(\Set), \Set] \to \Sh(\textrm{f.p.} {\mathbb T}\textrm{-mod}(\Set)^{\textrm{op}}, J)$ is the associated sheaf functor.\\
Given $\phi(x^{A_{1}}, \ldots, x^{A_{n}})$, suppose that $P_{\{\vec{x}. \phi\}}$ is definable by a geometric formula $\chi$ over $\Sigma$; then, by Theorem \ref{presheaf}, $[[\phi]]_{M_{\mathbb T}}$ and $[[\chi]]_{M_{\mathbb T}}$ are complemented subobjects in the lattice $\Sub_{\cal E}(M_{\mathbb T}A_{1}\times \ldots M_{\mathbb T}A_{n})$. Now, the fact that $a_{J}$ is a geometric functor implies that it preserves the interpretation of all geometric formulae, from which it follows that $[[\phi]]_{U'}$ and $[[\chi]]_{U'}$ are complemented subobjects in the lattice $\Sub_{\cal F}(U'A_{1}\times \cdots \times U'A_{n})$; then, from the conservativity of $U'$ and Lemma \ref{lemmafond}, it follows that $\phi \vee \neg^{{\mathbb T}'} \phi \stackrel{{\mathbb T}'}{\sim} \top$, as required. This proves the thesis when $P$ is the property $P_{\{\vec{x}. \phi\}}$. To prove the theorem for the property $Z^{{\mathbb T}'}_{\{\vec{x}.\phi\}}$, it suffices to invoke the thesis for the property $P_{\{\vec{x}. \phi \vee \neg^{{\mathbb T}'}\phi\}}$ and observe that in any Heyting algebra $H$, in particular in those of the form $\Sub_{{\cal C}_{{\mathbb T}'}}(\{\vec{x}. \top\})$, for any element $a\in H$, $a\vee \neg a=1$ if and only if $(a \vee \neg a)\vee \neg (a\vee \neg a)=1$.\\
(ii) It is immediate from Theorem \ref{presheaf} that if $\phi \vee \neg^{{\mathbb T}} \phi \stackrel{\mathbb T}{\sim} \top$ then $P_{\{\vec{x}. \phi\}}$ is definable by $\neg \phi$.     
\end{proofs}

\begin{corollary}
Let $\mathbb T$ be a theory of presheaf type over a signature $\Sigma$. Then\\
(i) the category $\textrm{f.p.} {\mathbb T}\textrm{-mod}(\Set)$ is a groupoid if and only if for every geometric formula $\phi$ over $\Sigma$, $P_{\{\vec{x}. \phi\}}$ is definable by a geometric formula over $\Sigma$;\\
(ii) the category $\textrm{f.p.} {\mathbb T}\textrm{-mod}(\Set)^{\textrm{op}}$ satisfies the right Ore condition if and only if for every geometric formula $\phi$ over $\Sigma$, $Q^{\mathbb T}_{\{\vec{x}. \phi\}}$ is definable by a geometric formula over $\Sigma$.
\end{corollary} 

\begin{proofs}
We have the following chain of equivalences: for every geometric formula $\phi$ over $\Sigma$, $P_{\{\vec{x}. \phi\}}$ (resp. $Q^{\mathbb T}_{\{\vec{x}. \phi\}}$) is definable by a geometric formula over $\Sigma$ if and only if (by Theorem \ref{teofond}) $\phi \vee \neg^{\mathbb T} \phi \stackrel{\mathbb T}{\sim} \top$ (resp. $\neg^{\mathbb T} \phi \vee \neg^{\mathbb T} \neg^{\mathbb T} \phi \stackrel{\mathbb T}{\sim} \top$), if and only if (cfr. \cite{OC3}) $[\textrm{f.p.} {\mathbb T}\textrm{-mod}(\Set), \Set]$ is Boolean (resp. De Morgan), if and only if (cfr. \cite{OC3}) $\textrm{f.p.} {\mathbb T}\textrm{-mod}(\Set)^{\textrm{op}}$ is a groupoid (resp. satisfies the right Ore condition). 
\end{proofs}

\section{Booleanizations and DeMorganizations via universal models}

In \cite{OC3} we provided essentially two approaches to investigate whether the classifying topos of a geometric theory $\mathbb T$ is Boolean (resp. De Morgan) or not, one of logical nature based on a syntactic characterization of the relevant class of geometric theories and the other one applicable in connection with a given site of definition of the classifying topos; accordingly, the process of Booleanization (resp. DeMorganization) was described both at the level of theories and in terms of the given site of definition of the topos. The purpose of this section is to relate these two approaches to one another by adopting a third point of view, that of universal models.\\  
Specifically, the syntactic criteria established in \cite{OC3} assert that a geometric theory $\mathbb T$ over a signature $\Sigma$ is classified by a Boolean (resp. De Morgan) topos if and only if for any geometric formula $\phi$ over $\Sigma$, $\phi \vee \neg^{\mathbb T} \phi \stackrel{{\mathbb T}}{\sim} \top$ (resp. $\neg^{\mathbb T} \phi \vee \neg^{\mathbb T}\neg^{\mathbb T} \phi \stackrel{{\mathbb T}}{\sim} \top$).\\
We recall that the dense topology on a category $\cal C$ is the Grothendieck topology $D_{\cal C}$ on $\cal C$ whose $D_{\cal C}$-covering sieves are exactly the stably non-empty ones (equivalently, the sieves of the form $R \vee \neg R$ for a sieve $R$ in $\cal C$) and the De Morgan topology on $\cal C$ defined in \cite{OC3} is the Grothendieck topology $M_{\cal C}$ on $\cal C$ generated by the sieves of the form $\neg R \vee \neg \neg R$ for a sieve $R$ in $\cal C$.\\  
The general criterion for a Grothendieck topos $\Sh({\cal C}, J)$ to be Boolean (resp. De Morgan) that we proved in \cite{OC3} is expressed in terms of the `reduced site' $(\tilde{\cal C}, J|_{\tilde{\cal C}})$, where $\tilde{\cal C}$ is the full subcategory of $\cal C$ on the objects which are not $J$-covered by the empty sieve, and asserts that $\Sh({\cal C}, J)$ is Boolean (resp. De Morgan) if and only if $J|_{\tilde{\cal C}}=D_{\tilde{\cal C}}$ (resp. $M_{\tilde{\cal C}}\leq J|_{\tilde{\cal C}}$).\\   

The following proposition describes the effect of taking subtoposes on the general and the reduced site representations of a given Grothendieck topos.\\
\begin{proposition}
Let $\cal C$ be a category and $\cal D$ a subcategory of $\cal C$. For a Grothendieck topology $J'$ on $\cal C$, denote by $J'|_{\cal D}$ the induced Grothendieck topology on $\cal D$, and for a Grothendieck topology $Z$ on $\cal D$, denote by $\overline{Z}$ the Grothendieck topology on $\cal C$ defined by:\\
$S\in \overline{Z}(c)$ if and only if for any arrow $f:d\to c$ in $\cal C$ with $d\in {\cal D}$, $f^{\ast}(C)\cap arr({\cal D})\in Z(c)$. Then:\\ 
(i) If $\cal D$ is $J$-dense then the assigments $J'\to J'|_{\cal D}$ and $Z\to \overline{Z}$ define a bijection between the class of Grothendieck topologies on $\cal C$ which contain $J$ and the class of Grothendieck topologies on $\cal D$ which contain $J|_{\cal D}$;\\
(ii) If $\cal C$ is locally small and $\cal D$ is small then for each Grothendieck topology $J'$ on $\cal C$ such that $J'\supseteq J$, the Comparison Lemma yields an equivalence $\Sh({\cal C}, J')\simeq \Sh({\cal D}, J'|_{\cal D})$.   
\end{proposition} 

\begin{proofs}
First, observe that if $J'\supseteq J$ and $\cal D$ is $J$-dense then $\cal D$ is also $J'$-dense; from this, part (ii) of the proposition immediately follows. So it remains to prove the following two facts:\\
(i) if $J$ is a Grothendieck topology on $\cal C$ and $\cal D$ is a $J$-dense subcategory of $\cal C$ then for any sieve $S$ in $\cal C$ on an object $c\in {\cal C}$, $S\in J(c)$ if and only if for any arrow $f:d\to c$ in $\cal C$ with $d\in {\cal D}$, $f^{\ast}(S)\cap arr({\cal D})\in J|_{\cal D}(c)$;\\
(ii) if $Z$ is a Grothendieck topology on $\cal D$ then $Z=\overline{Z}|_{\cal D}$.\\ 
Fact (i) was already observed in \cite{El2}, but we give a proof for the reader's convenience. One direction is obvious; let us then prove the other one. If for any arrow $f:d\to c$ in $\cal C$ with $d\in {\cal D}$, $f^{\ast}(S)\cap arr({\cal D})\in J|_{\cal D}(c)$ then by Lemma C2.2.2(i) \cite{El2} $f^{\ast}(S)$ is $J$-covering; hence the thesis follows from the transitivity axiom for Grothendieck topologies by invoking property (i) in the definition of $J$-dense subcategory given p. 546 \cite{El2}.\\
To prove fact (ii), we recall that for any sieve $R$ in $\cal D$ on $c\in {\cal D}$, the sieve $\overline{R}$ generated by $R$ in $\cal C$ is $\overline{Z}$-covering if and only if  $R\in \overline{Z}|_{\cal D}(c)$ (Lemma C2.2.2(i) \cite{El2}). Thus we have that $R\in \overline{Z}|_{\cal D}(c)$ if and only if $\overline{R}\in \overline{Z}(c)$, if and only if for any arrow $f:d\to c$ in $\cal C$ with $d\in {\cal D}$ $f^{\ast}(\overline{R})\cap arr({\cal D})\in Z(d)$, if and only if (take $f=1_{c}$) $R=\overline{R}\cap arr({\cal D})\in Z(c)$.   
\end{proofs}
A notable application of the proposition is in the following context: given a small category $\cal C$ and a Grothendieck topology $J$ on $\cal C$, the full subcategory $\tilde{\cal C}$ of $\cal C$ on the objects which are not $J$-covered by the empty sieve is $J$-dense and hence the subtoposes of $\Sh({\cal C}, J)$ can be identified with the subtoposes of $\Sh(\tilde{\cal C}, J|_{\tilde{\cal C}})$. Now, let us recall from \cite{OC3} that the Booleanization (resp. DeMorganization) of an elementary topos $\cal E$ is the largest dense Boolean subtopos of $\cal E$ (resp. the largest dense subtopos of $\cal E$ satisfying De Morgan's law); the corresponding local operator on $\cal E$ is the double-negation topology on $\cal E$ (resp. the De Morgan topology on $\cal E$). In fact, as we saw in \cite{OC3}, the Booleanization and DeMorganization of a topos $\Sh({\cal C}, J)$ are more naturally described in terms of the representation $\Sh(\tilde{\cal C}, J|_{\tilde{\cal C}})$. Specifically, by Proposition 2.5 \cite{OC3} and its well-known analogue for the Booleanization, the Booleanization (resp. DeMorganization) of $\Sh(\tilde{\cal C}, J|_{\tilde{\cal C}})$ is the topos $\Sh(\tilde{\cal C}, J|_{\tilde{\cal C}}\vee M_{\tilde{\cal C}})$ (resp. $\Sh(\tilde{\cal C}, D_{\tilde{\cal C}})$).\\   
From the proposition above, we can then deduce that the Booleanization of $\Sh({\cal C}, J)$ is the topos $\Sh({\cal C}, J_{b})$ where $J_{b}$ is the Grothendieck topology on $\cal C$ defined by: $\emptyset \in J_{b}(c)$ for $c\notin \tilde{\cal C}$, $S\in J_{b}(c)$ if and only if $S\supseteq \overline{T}$ for some $T\in D_{\tilde{{\cal C}}}(c)$, for $c\in \tilde{\cal C}$; similarly, the DeMorganization of $\Sh({\cal C}, J)$ is the topos $\Sh({\cal C}, J_{m})$ where $J_{m}$ is the Grothendieck topology on $\cal C$ defined by: $\emptyset \in J_{m}(c)$ for $c\notin \tilde{\cal C}$, $S\in J_{m}(c)$ if and only if $S\supseteq \overline{T}$ for some $T\in (M_{\tilde{{\cal C}}}\vee J)(c)$, for $c\in \tilde{\cal C}$.\\

Now, let us present an alternative view of Booleanizations and DeMorganizations.

\begin{proposition}\label{elem}
Let $\cal E$ be an elementary topos. Then:\\
(i) The double negation topology $\neg\neg$ on $\cal E$ is the smallest local operator $b$ on $\cal E$ such that all the monomorphisms of the form $A\vee \neg A\mono E$ for a subobject $A\mono E$ in $\cal E$ are $b$-dense;\\
(ii) The De Morgan topology on $\cal E$ is the smallest local operator $m$ on $\cal E$ such that all the monomorphisms of the form $\neg A\vee \neg \neg A\mono E$ for a subobject $A\mono E$ in $\cal E$ are $m$-dense.
\end{proposition}

\begin{proofs}
(i) By Theorem 2.4 \cite{OC3}, the double negation topology $\neg \neg$ on $\cal E$ is the smallest local operator $j$ on $\cal E$ such that the canonical monomorphism $(\top, \bot):2=1\amalg 1\rightarrowtail \Omega$ is $j$-dense. Now, observe that for any subobject $A\mono E$ in $\cal E$, the subobject $ A\vee \neg A\mono E$ is the pullback of $(\top, \bot)$ along the classifying map of $A\mono E$; so if $(\top, \bot)$ is $j$-dense for a local operator $j$ on $\cal E$ then also $A\vee \neg A\mono E$ is $j$-dense (recall that a monomorphism $l$ is $j$-dense if and only if the associated sheaf functor $a_{j}$ sends $l$ to an isomorphim, so that our claim follows from the fact that associated sheaf functors preserve pullbacks). Conversely, note that the subobject $(\top, \bot)$ is itself of the form $A\vee \neg A\mono E$, for $A\mono E$ equal to $\top:1\mono \Omega$.\\ 

(ii) The De Morgan topology $m$ on $\cal E$ is by definition (see \cite{OC3}) the smallest local operator $j$ on $\cal E$ such that the canonical monomorphism $(\top, \bot):2=1\amalg 1\rightarrowtail \Omega_{\neg \neg}$ is $j$-dense. Let us argue similarly as above. For any subobject $A\mono E$ in $\cal E$, the subobject $\neg A\vee \neg \neg A\mono E$ is the pullback of $(\top, \bot)$ along the classifying map of $\neg A\mono E$; so if $(\top, \bot)$ is $j$-dense for a local operator $j$ on $\cal E$ then also $\neg A\vee \neg \neg A\mono E$ is $j$-dense. Conversely, note that the subobject $(\top, \bot)$ is itself of the form $\neg A\vee \neg \neg A\mono E$, for $A\mono E$ equal to $\top:1\mono \Omega_{\neg \neg}$.   
\end{proofs}

In order to highlight the relationship between the definition of double negation topology (resp. De Morgan topology) on an elementary topos given in \cite{OC3} and the characterization of Proposition \ref{elem}, we observe a couple of useful facts (we acknowledge that the following result was stated without proof in \cite{PT}).\\

\begin{proposition}\label{negaz}
Let $\cal E$ be an elementary topos and $j\leq \neg\neg$ a local operator on $\cal E$ with associated sheaf functor $a_{j}:{\cal E}\to \sh_{j}({\cal E})$. Let $\neg$ and $\neg^{j}$ denote respectively the operation of pseudocomplementation of subobjects in $\cal E$ and in $\sh_{j}({\cal E})$. Then for any subobject $A\mono E$ in $\cal E$, $a_{j}(\neg A)\cong \neg^{j} a_{j}(A)$. In particular, if $A\mono E$ lies in $\sh_{j}({\cal E})$ then $\neg^{j} A \cong \neg A$. 
\end{proposition}

\begin{proofs}
The thesis follows as a consequence of the following two facts.\\
(1) It is well-known that $j\leq \neg\neg$ if and only if the associated sheaf functor $a_{j}$ sends non-zero subobjects to non-zero subobjects.\\
(2) Every subobject $D\mono a_{j}(E)$ in $\sh_{j}({\cal E})$ of an object of the form $a_{j}(E)$ is of the form $a_{j}(C)$ for some subobject $C\mono E$ in $\cal E$ (indeed, $D\mono a_{j}(E)$ is $j$-closed as a subobject in $\cal E$ and hence its classifying map factors through $\Omega_{j}\mono \Omega$; if $\chi:a_{j}(E)\to \Omega_{j}$ is this factorization then the composite of the map $E\to \Omega_{j}$ corresponding to $\chi$ via the adjunction with the monomorphism $\Omega_{j}\mono \Omega$ is easily seen to classify a subobject $C\mono E$ satisfying $a_{j}(C)\cong D$).\\
Now, the inequality $a_{j}(\neg A)\leq \neg^{j} a_{j}(A)$ holds for any $j$ (since $a_{j}$ preserves the zero-object and intersections of subobjects) while the inequality $\neg^{j} a_{j}(A) \leq a_{j}(\neg A)$ can be proved as follows. By (2), the subobject $\neg^{j} a_{j}(A)\mono a_{j}(E)$ is of the form $a_{j}(C)$ for some $C\mono E$; so we have $a_{j}(C\cap A)=a_{j}(C)\cap a_{j}(A)=\neg^{j} a_{j}(A) \cap a_{j}(A)=0$ and hence by $(1)$ $C\cap A=0$ i.e. $C\leq \neg A$ from which it follows that $a_{j}(C)\leq a_{j}(\neg A)$, as required.\\
The last part of the proposition follows from the first by observing that if $A\mono E$ lies in $\sh_{j}({\cal E})$ then $A \cong a_{j}(A)$ and hence $\neg^{j} A \cong a_{j}(\neg A)$; but $\neg A\mono E$ is $\neg\neg$-closed and hence $j$-closed (since $j\leq \neg\neg$), from which it follows that $a_{j}(\neg A) \cong \neg A$.        
\end{proofs}

In particular, the proposition shows that the pseudocomplementation of a subobject in a Grothendieck topos $\Sh({\cal C}, J)$ corresponds via the equivalence $\Sh({\cal C}, J)\simeq \Sh(\tilde{\cal C}, J|_{\tilde{\cal C}})$ to the pseudocomplementation in the presheaf topos $[\tilde{\cal C}^{\textrm{op}}, \Set]$ of the image of it under the equivalence.\\  

\begin{proposition}\label{rev}
Let $\cal E$ be a Grothendieck topos $\Sh({\cal C}, J)$ and $m:A\mono E$ a monomorphism in $\cal E$. Given a local operator $j'$ on $\cal E$ corresponding to a Grothendieck topology $J'\supseteq J$ on $\cal C$, $A\mono E$ is $j'$-dense if and only if for each $e\in E(c)$ the sieve $S^{A}_{(c, e)}:=\{f:d\to c \textrm{ | } E(f)(e)\in A(d)\}$ in $\cal C$ is $J'$-covering. In particular, the smallest local operator $j'$ on $\cal E$ such that $A\mono E$ is $j'$-dense corresponds to the Grothendieck topology $J'\supseteq J$ generated over $J$ by the collection of sieves $S^{A}_{(c,e)}$ as $(c,e)$ varies in the category of elements of the functor $E$.   
\end{proposition}  

\begin{proofs}
Let us denote by $a_{j'}:{\cal E}\to \sh_{j'}({\cal E})$, $a_{J}:[{\cal C}^{\textrm{op}}, \Set]\to \Sh({\cal C}, J)$ and $a_{J'}:[{\cal C}^{\textrm{op}}, \Set]\to \Sh({\cal C}, J')$ the associated sheaf functors; then $m$ is $j'$-dense if and only if $a_{j'}(m)$ is an isomorphism. But $a_{j'}(m)\cong a_{j'}(a_{J}(m))\cong a_{J'}(m)$, where $m$ is regarded here as a subobject in $[{\cal C}^{\textrm{op}}, \Set]$. Now, $a_{J'}(m)$ is an isomorphism if and only if $m$ is $c_{J'}$-dense, where $c_{J'}$ is the universal closure operator on $[{\cal C}^{\textrm{op}}, \Set]$ corresponding to $J'$, if and only if $c_{J'}(m)=E$ i.e. for any $e\in E(c)$, $S^{A}_{(c, e)}\in J'(c)$ (cfr. formula (6) p. 235 \cite{MM}).
\end{proofs}

\begin{rmk}
\emph{Or course, the proposition assumes a particularly simple form if $E$ is representable, say $E={\cal C}(-,c)$; in fact, the monomorphisms $A\mono E$ can be identified with sieves $S$ on $c$ and hence the smallest local operator $j'$ on $\cal E$ such that $S\mono E$ is $j'$-dense corresponds to the Grothendieck topology $J'\supseteq J$ generated over $J$ by the sieve $S$. It is useful to keep this in mind while working with universal models; for example, if $\mathbb T$ is cartesian then $M_{\mathbb T}A_{1}\times \cdots M_{\mathbb T}A_{n}$ is representable (for any sorts $A_{1}, \ldots, A_{n}$).}     
\end{rmk}

By the results in \cite{OC3}, given an elementary topos $\cal E$ and a local operator $j$ on $\cal E$ such that $j\leq \neg\neg$, the Booleanization (resp. DeMorganization) of the topos $\sh_{j}(\cal E)$ is given by the subtopos $\sh_{\neg\neg}({\cal E})$ (resp. $\sh_{j\vee m}({\cal E})$), where $m$ is the De Morgan topology on $\cal E$. Now, these descriptions are easily seen, by using Proposition \ref{negaz}, to be equivalent to the characterizations given in Proposition \ref{elem}.\\
In the case $\cal E$ is a Grothendieck topos $\Sh({\cal C}, J)$, Proposition \ref{rev} provides another link between the two approaches; indeed, by the proof of Proposition \ref{elem}, the Booleanization (resp. DeMorganization) of $\Sh({\cal C}, J)$ corresponds to the the smallest local operator $l$ on $\cal E$ such that the monomorphism $2 \mono \Omega$ (resp. $2 \mono \Omega_{\neg\neg}$) in $\cal E$ is $l$-dense, and if $J\leq \neg\neg$ then $S^{2}_{(c, R)}=R\vee \neg R$ (resp. $S^{2}_{(c, R)}=\neg R\vee \neg \neg R$).\\

Let us recall from \cite{OC3} that, given a Grothendieck topos $\cal E$ with a separating set $\cal G$, $\cal E$ is Boolean (resp. De Morgan) if and only if all the subobject lattices of the form $\Sub_{\cal E}(c)$ for $c\in {\cal G}$ are Boolean (resp. De Morgan) algebras. This fact leads to the following result.\\

\begin{proposition}\label{gen}
Let $\cal E$ be a Grothendieck topos with a separating set $\cal G$. Then:\\
(i) The double negation topology $\neg\neg$ on $\cal E$ is the smallest local operator $b$ on $\cal E$ such that all the monomorphisms of the form $A\vee \neg A\mono E$ for a subobject $A\mono E$ in $\cal E$ with $E\in {\cal G}$ are $b$-dense;\\
(ii) The De Morgan topology on $\cal E$ is the smallest local operator $m$ on $\cal E$ such that all the monomorphisms of the form $\neg A\vee \neg \neg A\mono E$ for a subobject $A\mono E$ in $\cal E$ with $E\in {\cal G}$ are $m$-dense.
\end{proposition}

\begin{proofs}
Let $b$ be the smallest local operator $j$ on $\cal E$ such that all the monomorphisms of the form $A\vee \neg A\mono E$ for a subobject $A\mono E$ in $\cal E$ with $E\in {\cal G}$ are $j$-dense. Then $b\leq \neg\neg$ and, by Proposition \ref{negaz} and the criterion stated before the proposition, the topos $\sh_{b}({\cal E})$ is Boolean; thus, since $\neg\neg$ is the smallest dense local operator $j$ on $\cal E$ such that $\sh_{j}({\cal E})$ is Boolean, $b=\neg\neg$, as required.\\ 
Part (ii) follows analogously to Part (i), by recalling that the De Morgan topology on $\cal E$ is the smallest dense local operator $j$ on $\cal E$ such that $\sh_{j}({\cal E})$ satisfies De Morgan's law.
\end{proofs}

Proposition \ref{gen}, together with Proposition \ref{separ}, provides us with a convenient approach to test whether the classifying topos $\cal E$ of a given geometric theory $\mathbb T$ is Boolean (resp. De Morgan) or not: ${\cal E}$ is Boolean (resp. De Morgan) if and only if all the subobject lattices of the form $\Sub_{\cal E}(MA_{1}\times \cdots \times MA_{n})$ where $M$ is a universal model of $\mathbb T$ in $\cal E$ are Boolean (resp. De Morgan) algebras.\\ 
This approach is also conceptually enlightening since, by Theorem \ref{definab}, universal models are unifying elements of syntax and semantics. More precisely, the link between it and our syntactic criteria for a geometric theory to be Boolean (resp. De Morgan) is exactly provided by Theorem \ref{unif}. In fact, forcing a property of subobjects in $\cal E$ of $MA_{1}\times \cdots \times MA_{n}$ to hold in a subtopos $\sh_{j}({\cal E})$ of $\cal E$ exactly corresponds to forcing the logical property of geometric formulae corresponding to them via Theorem \ref{unif} to hold in the quotient ${\mathbb T}'$ of $\mathbb T$ corresponding to the subtopos $\sh_{j}({\cal E})$ via Theorem 3.6 \cite{OC6}; indeed, by Lemma \ref{subtopos}, the image of $M$ under the associated sheaf functor $a_{j}$ is a universal, and hence conservative, model of ${\mathbb T}'$. For example, forcing $B \vee \neg B \mono MA_{1}\times \cdots \times MA_{n}$ to be an isomorphism corresponds (by Lemma \ref{lemmafond}) to forcing $\phi \vee \neg^{\mathbb T} \phi$ to be ${\mathbb T}'$-provably equivalent to $\top$, where $B\mono MA_{1}\times \cdots \times MA_{n}$ is the interpretation of $\phi$ in $M$.\\

For another application of Proposition \ref{gen}, consider a quotient ${\mathbb T}'$ of a theory of presheaf type $\mathbb T$ such that its associated $\mathbb T$-topology $J$ on $\textrm{f.p.} {\mathbb T}\textrm{-mod}(\Set)^{\textrm{op}}$ is subcanonical. We can describe the Booleanization (resp. DeMorganization) of the classifying topos $\Sh(\textrm{f.p.} {\mathbb T}\textrm{-mod}(\Set)^{\textrm{op}}, J)$ of ${\mathbb T}'$ as follows. $J_{b}$ (resp. $J_{m}$) is generated over $J$ by the sieves of the form $S\vee \neg S$ (resp. $\neg S \vee \neg \neg S$) for a $J$-closed sieve $S$ in $\textrm{f.p.} {\mathbb T}\textrm{-mod}(\Set)^{\textrm{op}}$, where the lattice operations $\vee$ and $\neg$ on subobjects are taken in $\Sh(\textrm{f.p.} {\mathbb T}\textrm{-mod}(\Set)^{\textrm{op}}, J)$.\\ 

Now our aim is to relate the approach to Booleanizations (resp. DeMorganizations) of Grothendieck toposes based on the notion of dense topology (resp. De Morgan topology) on a category to the point of view of universal models presented in this paper. Specifically, we show how it is possible, under appropriate hypotheses, to naturally convert a counterexample to a Grothendieck topos $\Sh({\cal C}, J)$ being Boolean (resp. De Morgan) i.e. a sieve $R$ in $\cal C$ such that $R \vee \neg R\notin J|_{\tilde{{\cal C}}}$ (resp. $\neg R \vee \neg \neg R \notin J|_{\tilde{{\cal C}}}$) into a logical counterexample to a theory $\mathbb T$ classified by $\Sh({\cal C}, J)$ being Boolean (resp. De Morgan) i.e. a geometric formula $\phi$ over the signature of $\mathbb T$ such that $\phi \vee \neg^{\mathbb T} \phi \stackrel{\mathbb T}{\not \sim} \top$ (resp. $\neg^{\mathbb T}\phi \vee \neg^{\mathbb T} \neg^{\mathbb T} \phi \stackrel{\mathbb T}{\not \sim} \top$).\\

Let $({\cal C}, J)$ be a Grothendieck site and $\tilde{\cal C}$ the full subcategory of $\cal C$ on the objects which are not $J$-covered by the empty sieve. Recall that, by the Comparison Lemma, $\Sh({\cal C}, J)\simeq \Sh(\tilde{\cal C}, J|_{\tilde{\cal C}})$. As observed in Example C2.2.4 \cite{El2}, the full embedding $\tilde{\cal C}^{\textrm{op}} \hookrightarrow {\cal C}^{\textrm{op}}$ gives rise to a geometric inclusion $[\tilde{\cal C}^{\textrm{op}}, \Set] \to [{\cal C}^{\textrm{op}}, \Set]$ whose inverse image is the obvious restriction map; this inclusion is isomorphic (via the Comparison Lemma) to the canonical inclusion $\Sh({\cal C}, J_{\cal C}) \hookrightarrow [{\cal C}^{\textrm{op}}, \Set]$ where $J_{\cal C}$ is the Grothendieck topology on $\cal C$ defined by saying that a sieve $R$ on $c$ is $J_{\cal C}$-covering if and only if it contains all the morphisms from objects of $\tilde{{\cal C}}$ to $c$.\\
Now, suppose that $\cal C$ is the opposite of the category of finitely presentable models of a theory of presheaf type $\mathbb T$ over a signature $\Sigma$. Then the subtopos $[\tilde{\cal C}^{\textrm{op}}, \Set] \simeq \Sh({\cal C}, J_{\cal C}) \hookrightarrow [{\cal C}^{\textrm{op}}, \Set]$ of the classifying topos $[{\cal C}^{\textrm{op}}, \Set]$ of $\mathbb T$ corresponds via Theorem 3.6 \cite{OC6} to a quotient ${\mathbb T}_{\tilde{\cal C}}$ of $\mathbb T$ classified by $[\tilde{\cal C}^{\textrm{op}}, \Set]$; we will refer to this quotient as the $\tilde{\cal C}$-reduct of $\mathbb T$. By Lemma \ref{subtopos}, the image $M_{{\mathbb T}_{\tilde{\cal C}}}:=a_{J_{\cal C}}(M_{\mathbb T})$ of the universal model $M_{\mathbb T}$ of $\mathbb T$ under the associated sheaf functor $a_{J_{\cal C}}:[{\cal C}^{\textrm{op}}, \Set] \to \Sh({\cal C}, J_{\cal C})$ is a universal model of ${\mathbb T}_{\tilde{\cal C}}$; but for any sort $A$, $a_{J_{\cal C}}(M_{\mathbb T})A$ is the restriction of the functor $M_{\mathbb T}A$ to the category $\tilde{\cal C}$ and hence Theorem \ref{presheaf} yields the following explicit description of the universal model $M_{{\mathbb T}_{\tilde{\cal C}}}$ of ${\mathbb T}_{\tilde{\cal C}}$ in $[\tilde{\cal C}^{\textrm{op}}, \Set]$: for any sort $A$ of $\Sigma$ and any $M\in \tilde{\cal C}$, $(M_{{\mathbb T}_{\tilde{\cal C}}}A)M=MA$.\\       
Now, $J$, as a Grothendieck topology on $\textrm{f.p.} {\mathbb T}\textrm{-mod}(\Set)^{\textrm{op}}$, is the associated $\mathbb T$-topology of a quotient ${\mathbb T}'$ of $\mathbb T$. We note that that for any geometric formula $\phi$ over $\Sigma$, $\neg^{{\mathbb T}'} \phi \stackrel{{\mathbb T }'}{\sim} \neg^{{\mathbb T}_{\tilde{\cal C}}} \phi$. Indeed, by definition of $\tilde{\cal C}$, $J|_{\tilde{\cal C}}$ is dense and hence, by Proposition \ref{negaz}, the associated sheaf functor $a_{J|_{\tilde{\cal C}}}:[\tilde{{\cal C}}^{\textrm{op}}, \Set]\to \Sh(\tilde{{\cal C}}, J|_{\tilde{{\cal C}}})$ preserves pseudocomplementation of subobjects. But $M_{{\mathbb T}'}:=a_{J|_{\tilde{\cal C}}}(M_{{\mathbb T}_{\tilde{\cal C}}})$ is a universal model of ${\mathbb T}'$ and hence $[[\neg^{{\mathbb T}'} \phi]]_{M_{{\mathbb T}'}}=\neg [[\phi]]_{M_{{\mathbb T}'}}=\neg a_{J|_{\tilde{\cal C}}}([[\phi]]_{M_{{\mathbb T}_{\tilde{\cal C}}}})=a_{J|_{\tilde{\cal C}}}(\neg [[\phi]]_{M_{{\mathbb T}_{\tilde{\cal C}}}})=a_{J|_{\tilde{\cal C}}}([[\neg^{{\mathbb T}_{\tilde{\cal C}}}\phi]]_{M_{{\mathbb T}_{\tilde{\cal C}}}})=[[\neg^{{\mathbb T}_{\tilde{\cal C}}}\phi]]_{M_{{\mathbb T}'}}$, where the first and fourth equalities follow from Lemma \ref{lemmafond} and the second and fifth equalities follow from the fact that associated sheaf functors preserve the interpretation of geometric formulae. Therefore, by the conservativity of universal models, it follows that $\neg^{{\mathbb T}'} \phi \stackrel{{\mathbb T }'}{\sim} \neg^{{\mathbb T}_{\tilde{\cal C}}} \phi$, as required. We note that, passing from ${\mathbb T}_{\tilde{\cal C}}$ to its Booleanization (resp. DeMorganization) ${\mathbb T}''$ forces all the properties of the form $Z^{{\mathbb T}_{\tilde{\cal C}}}_{\phi}$ (resp. $Z^{{\mathbb T}_{\tilde{\cal C}}}_{\neg^{{\mathbb T}_{\tilde{\cal C}}} \phi}$) not to hold in any ${\mathbb T}''$-model in a Grothendieck topos.\\

Let us now suppose that $\mathbb T$ is cartesian and that all the representing objects $F_{A_{1}, \ldots, A_{n}}$ of the functors $M_{\mathbb T}A_{1}\times \cdots \times M_{\mathbb T}A_{n}:\textrm{f.p.} {\mathbb T}\textrm{-mod}(\Set) \to \Set$ (for any sorts $A_{1}, \ldots, A_{n}$ of $\Sigma$) lie in $\tilde{\cal C}$; under these hypotheses, we can easily turn a `categorical' counterexample to the classifying topos of ${\mathbb T}'$ being Boolean (resp. De Morgan) into a logical counterexample to ${\mathbb T}'$ being Boolean (resp. De Morgan).\\
First, we note that if the classifying topos of ${\mathbb T}'$ is not Boolean (resp. De Morgan) then there must be a sieve $R$ on an object of the form $F_{A_{1}, \ldots, A_{n}}$ providing a countexample, since any object in the (cartesian) syntactic category of $\mathbb T$ injects into an object of the form $\{x^{A_{1}}, \ldots, x^{A_{n}}. \top\}$ (pullback functors, being logical, preserve all the Heyting operations between subobjects, cfr. also Remark 1.6(b) \cite{OC3}).\\ 
Then, given a sieve $R$ in $\tilde{\cal C}$ on $F_{A_{1}, \ldots, A_{n}}$, by logical completeness of universal models, $R$ is the interpretation in the universal model of ${\mathbb T}_{\tilde{\cal C}}$ in $[\tilde{\cal C}^{\textrm{op}}, \Set]$ of a geometric formula $\phi(x^{A_{1}}, \ldots, x^{A_{n}})$ over $\Sigma$. Now, if $R \vee \neg R\notin J|_{\tilde{\cal C}}(F_{A_{1}, \ldots, A_{n}})$ (resp. $\neg R\vee \neg \neg R \notin J|_{\tilde{\cal C}}(F_{A_{1}, \ldots, A_{n}})$) then, by the preceding discussion, $\phi \vee \neg^{{\mathbb T}'} \phi \stackrel{{\mathbb T }'}{\not \sim} \top$ (resp. $\neg^{{\mathbb T}'}\phi \vee \neg^{{\mathbb T}'} \neg^{{\mathbb T}'} \phi \stackrel{{\mathbb T }'}{\not \sim} \top$).\\   
Another interesting case to consider is the theory $\mathbb T$ of flat functors on $\cal C$ and its quotient ${\mathbb T }'$ of $J$-continuous flat functors on $\cal C$, which, by Diaconescu's theorem, is classified by the topos $\Sh({\cal C}, J)$. The universal model $U'$ of ${\mathbb T }'$ in $\Sh({\cal C}, J)$ assigns to each sort $c\in {\cal C}$ of the signature of ${\mathbb T }'$ the functor $a_{J}({\cal C}(-, c))$, where $a_{J}:[{\cal C}^{\textrm{op}}, \Set]\to \Sh({\cal C}, J)$ is the associated sheaf functor. Now, given a sieve $R$ in $\tilde{{\cal C}}$ such that $R \vee \neg R$ is not $J|_{\tilde{{\cal C}}}$-covering (resp. $\neg R \vee \neg\neg R$ is not $J|_{\tilde{{\cal C}}}$-covering), we have by Proposition \ref{negaz} that $a_{J}(R) \vee \neg a_{J}(R)\cong a_{J}(R \vee \neg R)\ncong a_{J}({\cal C}(-, c))$ (resp. $\neg a_{J}(R) \vee \neg \neg a_{J}(R)\cong a_{J}(\neg R \vee \neg\neg R)\ncong a_{J}({\cal C}(-, c))$); hence the geometric formula $\phi$ over the signature of ${\mathbb T }'$ whose interpretation in the universal model $U'$ is $a_{J}(R)$ satisfies $\phi \vee \neg^{{\mathbb T}'} \phi \stackrel{{\mathbb T }'}{\not \sim} \top$ (resp. $\neg^{{\mathbb T}'}\phi \vee \neg^{{\mathbb T}'} \neg^{{\mathbb T}'} \phi \stackrel{{\mathbb T }'}{\not \sim} \top$), thus providing a counterexample to the theory ${\mathbb T }'$ being Boolean (resp. De Morgan).\\
  
Given a general Grothendieck site $({\cal C}, J)$, a notable advantage of working with the reduced site $(\tilde{\cal C}, J|_{\tilde{\cal C}})$ is that, since $J|_{\tilde{\cal C}}\leq \neg\neg_{\tilde{\cal C}}$, every sieve in $\tilde{\cal C}$ of the form $\neg R$ is automatically $J|_{\tilde{\cal C}}$-closed (being $\neg\neg$-closed); an application of this fact is given by the following proposition.

\begin{proposition}\label{chiusura}
Let $\mathbb T$ be a theory of presheaf type with universal model $M_{\mathbb T}$ in $[\textrm{f.p.} {\mathbb T}\textrm{-mod}(\Set), \Set]$ (as in Theorem \ref{presheaf}) and $J$ a Grothendieck topology on $\textrm{f.p.} {\mathbb T}\textrm{-mod}(\Set)^{\textrm{op}}$ such that every $J$-covering sieve is non-empty. Then, given $\phi=\phi(x^{A_{1}}, \cdots, x^{A_{n}})$ and a subobject $E\mono M_{\mathbb T}A_{1}\times \cdots M_{\mathbb T}A_{n}$ in $[\textrm{f.p.} {\mathbb T}\textrm{-mod}(\Set), \Set]$ such that $[[\phi]]_{M_{\mathbb T}}=\neg E$ in $[\textrm{f.p.} {\mathbb T}\textrm{-mod}(\Set), \Set]$, for any $M\in \textrm{f.p.} {\mathbb T}\textrm{-mod}(\Set)$ and $\vec{a}\in MA_{1}\times \cdots MA_{n}$, $\{f:M\to N \textrm{ in f.p.} {\mathbb T}\textrm{-mod}(\Set) \textrm{ | } f(\vec{a})\in [[\phi]]_{N}\}\in J(M)$ implies $\vec{a}\in [[\phi]]_{M}$.
\end{proposition}

\begin{proofs}
This follows at once from the fact that any subobject in $[\textrm{f.p.} {\mathbb T}\textrm{-mod}(\Set), \Set]$ of the form $\neg E$ is $D$-closed, where $D$ is the dense topology on $\textrm{f.p.} {\mathbb T}\textrm{-mod}(\Set)^{\textrm{op}}$, and our hypothesis that $J\leq D$, by recalling the expression for $\overline{[[\vec{x}. \phi]]_{M_{\mathbb T}}}^{J}$ derived in section \ref{terza} above. 
\end{proofs}

Note that, by using the explicit expression for $M_{\mathbb T}$ given in Theorem \ref{presheaf}, the condition  $[[\phi]]_{M_{\mathbb T}}=\neg E$ rewrites as follows: for any $M\in \textrm{f.p.} {\mathbb T}\textrm{-mod}(\Set)$,

\[
\begin{array}{ccl}
[[\vec{x}.\phi]]_{M} & = & \{(a_{1}, \cdots, a_{n})\in MA_{1}\times \cdots MA_{n} \textrm{ | for all } f:M\to N \\
& & \textrm{ in f.p.} {\mathbb T}\textrm{-mod}(\Set) \textrm{ | } f(a_{1}, \cdots a_{n})\notin E(N)\}.   
\end{array}
\]

As we shall see in the next section, this proposition can be profitably used to establish `closure results' for various mathematical properties.\\
The proposition is most notably applied in connection to reducts of theories $\mathbb T$ of presheaf type such that the functor $M_{\mathbb T}A_{1}\times \cdots M_{\mathbb T}A_{n}:\tilde{\cal C}^{\textrm{op}}\to \Set$ is representable by an object $F_{A_{1}, \ldots, A_{n}}$, so that subobjects of $M_{\mathbb T}A_{1}\times \cdots M_{\mathbb T}A_{n}$ in $[{\tilde{{\cal C}}}^{\textrm{op}}, \Set]$ correspond to sieves on $F_{A_{1}, \ldots, A_{n}}$ in $\tilde{\cal C}$. In fact, the operation of pseudocomplementation of sieves in a presheaf topos $[{\cal C}^{\textrm{op}}, \Set]$ is uninteresting when the category $\cal C$ has a lot of structure; for example, if $\cal C$ satisfies the right Ore condition (notice that this is always the case if $\cal C$ is opposite of the category of models of a cartesian theory) then, given a sieve $R$ on $c\in {\cal C}$, $\neg R=\emptyset$ if $R$ is non-empty and $\neg R=M_{c}$ if $R$ is empty.\\

\section{Examples}\label{exa}

\hspace{0.7 cm} (i) Let ${\mathbb T}$ be the algebraic theory of commutative rings with unit. The signature $\Sigma$ of $\mathbb T$ is one-sorted and consists of two binary function symbols $+$ and $\cdot$, one unary function symbol $-$ and two constants $0$ and $1$. Let ${\mathbb T}'$ be the theory of non-trivial rings i.e. the quotient of $\mathbb T$ obtained by adding the sequent $(0=1 \vdash_{[]} \bot)$. Let ${\cal C}:=\textrm{f.p.} {\mathbb T}\textrm{-mod}(\Set)$ be the category of finitely presented (equivalently, finitely generated) rings and $\tilde{\cal C}$ be the full subcategory of it on the non-zero rings.\\  
Consider the following geometric formulae over $\Sigma$: $\phi_{1}(x)\equiv (\exists y)(x\cdot y =1)$, which says that an element is invertible, $\phi_{2}(x)\equiv \mathbin{\mathop{\textrm{\huge $\vee$}}\limits_{n\in {\mathbb N}}} (x^{n}=0)$, which says that an element is nilpotent, and $\phi_{3}(x)\equiv (x=0)$. Adopting the notation of section \ref{sec}, let us denote by $R_{1}$ the sieve in $\tilde{\cal C}^{\textrm{op}}$ generated by the single arrow $d([\phi_{1}]):M_{\{x.\top\}}\to M_{\{x. \phi_{1}\}}$. We proved in \cite{OC3} that $\neg R_{1}=\{f:M_{\{x.\top\}}\to M \textrm{ in $\tilde{\cal C}$ | } f(\xi_{\top})\in [[\phi_{2}]]_{M} \}$ and $\neg \neg R_{1}=\{f:M_{\{x.\top\}}\to M \textrm{ in $\tilde{\cal C}$ | } f(\xi_{\top})\in [[\phi_{1}]]_{M} \}=R_{1}$, where the pseudocomplementation $\neg$ is taken in the presheaf topos $[\tilde{\cal C}, \Set]$; similarly, one can prove that the sieve $R_{3}$ in $\tilde{\cal C}$ generated by the single arrow $d([\phi_{3}])$ satisfies $\neg R_{3}=\{f:M_{\{x.\top\}}\to M \textrm{ in $\tilde{\cal C}$ | } f(\xi_{\top})\in [[\phi_{1}]]_{M} \}=R_{1}$. Now, from the results of section \ref{sec}, we have that $[[\phi_{1}]]_{M_{{\mathbb T}'}}=R_{1}$, $[[\phi_{2}]]_{M_{{\mathbb T}'}}=\{f:M_{\{x.\top\}}\to M \textrm{ in $\tilde{\cal C}$ | } f(\xi_{\top})\in [[\phi_{2}]]_{M} \}=\neg R_{1}$ and $[[\phi_{3}]]_{M_{{\mathbb T}'}}=R_{3}$, where $M_{{\mathbb T}'}$ is the universal model of ${\mathbb T}'$ in $[\tilde{\cal C}, \Set]$. Hence $\neg^{{\mathbb T}'} \phi_{2} \stackrel{{\mathbb T }'}{\sim} \neg^{{\mathbb T}'} \phi_{3} \stackrel{{\mathbb T }'}{\sim}  \phi_{1}$ and $\neg^{{\mathbb T}'} \phi_{1} \stackrel{{\mathbb T }'}{\sim} \phi_{2}$. Notice that it follows from these relations that the sequent $\top \vdash_{x} \phi_{3} \vee \phi_{1}$ is provable in the Booleanization of ${\mathbb T}'$ and that the sequent $(\top \vdash_{x} (\exists y)(x\cdot y =1) \vee (\mathbin{\mathop{\textrm{\huge $\vee$}}\limits_{n\in {\mathbb N}}} (x^{n}=0)))$ is provable in the DeMorganization of ${\mathbb T}'$. Also, it follows from Theorem \ref{teofond} that the property of an element of a (finitely generated) ring to be neither invertible nor nilpotent is not definable by a geometric formula over $\Sigma$.\\  
Now, let us consider the following quotients of ${\mathbb T}'$: the theory ${\mathbb T}_{l}$ (resp. ${\mathbb T}_{d}$) of local rings (resp. of integral domains) as axiomatized in \cite{OC3}, and the coherent theory ${\mathbb T}_{f}$ of fields i.e. the quotient of ${\mathbb T}'$ obtained by adding to it the sequent $(\top \vdash_{x} (x=0) \vee (\exists y) (x\cdot y=1))$. The ${\mathbb T}'$-topologies on $\tilde{\cal C}^{\textrm{op}}$ associated to ${\mathbb T}_{l}$, ${\mathbb T}_{d}$ and ${\mathbb T}_{f}$ will be denoted respectively by $J_{l}$, $J_{d}$ and $J_{f}$. These topologies have the following descriptions:\\
$J_{l}$ is the restriction to $\tilde{\cal C}$ of the Zariski topology on $\cal C$ and hence is defined by: $S\in J_{l}(A)$ if and only if $S$ contains a finite family $\{\xi_{s_{i}}:A\rightarrow A[{s_{i}}^{-1}] \textrm{ | } 1\leq i\leq n\}$ of canonical inclusions $\xi_{s_{i}}:A\rightarrow A[{s_{i}}^{-1}]$ in $\textrm{f.p.} {\mathbb T}\textrm{-mod}(\Set)$ where $\{s_{1},\ldots, s_{n}\}$ is any set of non-nilpotent elements of $A$ which is not contained in any proper ideal of $A$ (cfr. \cite{OC3});\\
$J_{d}$ is given by: $S\in J_{d}(A)$ if and only if $S$ contains a finite family $\{\pi_{a_{i}}:A\rightarrow A/(a_{i}) \textrm{ | } 1\leq i\leq n\}$ of canonical projections $\pi_{a_{i}}:A\rightarrow A/(a_{i})$ in $\textrm{f.p.} {\mathbb T}\textrm{-mod}(\Set)$ where $\{a_{1},\ldots,a_{n}\}$ is any set of non-invertible elements of $A$ such that $a_{1}\cdot \ldots \cdot a_{n}=0$;\\
$J_{f}$ is generated by the sieves which contain families of the form $(\xi_{a}:A\to A[a^{-1}], \pi_{a}:A\rightarrow A/(a))$, where $a$ is a neither nilpotent nor invertible element of $A$.\\
Notice that all the Grothendieck topologies $J_{l}$, $J_{d}$ and $J_{f}$ are dense on $\tilde{\cal C}^{\textrm{op}}$; in particular, the Booleanizations of the theories ${\mathbb T}'$, ${\mathbb T}_{l}$, ${\mathbb T}_{d}$ all coincide with the Booleanization of the theory ${\mathbb T}_{f}$ of fields.\\
Let us now apply Proposition \ref{chiusura} in the context of our topologies $J_{l}$, $J_{d}$ and $J_{f}$; observe that, by the calculations above, one can take the property $\phi$ to be either $\phi_{1}$ or $\phi_{2}$.\\
For example, if $\phi$ is $\phi_{1}$ (resp. $\phi_{2}$) and $J$ is $J_{l}$ then we get the following result:\\ 
for any finitely generated ring $A$ and set $\{s_{1},\ldots, s_{n}\}$ of non-nilpotent elements of $A$ which is not contained in any proper ideal of $A$, for any $a\in A$, $a$ is invertible (resp. nilpotent) in $A$ if and only if for each $i=1, \ldots, n$, $\xi_{s_{i}}(a)$ is invertible (resp. nilpotent) in $A[{s_{i}}^{-1}]$.\\
If we take $\phi$ to be $\phi_{1}$ (resp. $\phi_{2}$) and $J$ to be $J_{f}$ then Proposition \ref{chiusura} gives:\\
for any finitely generated ring $A$ and any neither nilpotent nor invertible element $s\in A$, for any $a\in A$, $a$ is invertible (resp. nilpotent) in $A$ if and only if $\xi_{s}(a)$ is invertible (resp. nilpotent) in $A[{s}^{-1}]$ and $\pi_{s}(a)$ is invertible (resp. nilpotent) in $A/(s)$.\\
Let us now extract from the ideas underlying Proposition \ref{chiusura} a direct proof of the `invertible' version of the latter result.\\  
Given an element $a\in A$ such that $\xi_{s}(a)$ is invertible in $A[{s}^{-1}]$ and $\pi_{s}(a)$ is invertible in $A/(s)$, consider the ring $B:=A/(a)$ and the natural projection map $l:A\to B$. From the fact that we have commutative squares

\[  
\xymatrix {
A \ar[d]^{\pi_{s}} \ar[r]^{l}  & B \ar[d] & A \ar[d]^{\xi_{s}} \ar[r]^{l} &  B \ar[d] \\
A/(s) \ar[r] & B/(l(s)) & A[{s}^{-1}] \ar[r] & B[{l(s)}^{-1}]}
\]

in the category of (finitely generated) rings, it follows that both $B/(l(s))$ and $B[{l(s)}^{-1}]$ must be the zero-ring; indeed, our hypotheses imply that the image of $a$ under the diagonal arrow $A\to B/(l(s))$ (resp. $A\to B[{l(s)}^{-1}]$) is both an invertible and zero element of the ring $B/(l(s))$ (resp. $B[{l(s)}^{-1}]$). Then $B$ must itself be the zero ring i.e. $a$ must be inverible in $A$, as required.\\ 

Of course, each of the above-mentioned results can also be proved `concretely' by using algebraic manipulation, but we emphasize that the point of view offered by Proposition \ref{chiusura} has the advantage of providing a unifying framework for all these problems, which might otherwise  seem to be unrelated to each other.

Finally, notice that if $\mathbb Z$ is a quotient of ${\mathbb T}'$ having enough models and such that the associated ${\mathbb T}'$-topology of $\mathbb Z$ is dense on $\tilde{\cal C}$ (for example $\mathbb Z$ can be any of the quotients of ${\mathbb T}'$ which we considered above) then the syntactic relation $\stackrel{{\mathbb Z}}{\sim}$ between formulae can be rephrased semantically in terms of definable properties. For example, the fact that $\neg^{{\mathbb Z}} \phi_{3}$ ${\mathbb Z}$-provably implies $\phi_{1}$ rephrases as follows. Let $P(x)$ be a property of elements of $\mathbb Z$-models which is never satisfied by zero; if $P$ is definable by a geometric formula over the signature of rings then any element of a $\mathbb Z$-model satisfying $P$ is invertible.\\

\hspace{0.7 cm} (ii) Let $\mathbb T$ be the cartesian theory $\mathbb T$ of von Neumann regular rings. There is a notion of characteristic $char(R)$ of such a ring $R$; specifically,
\[
char(R)=\{p\in {\mathbb P}\cup \{0\} \textrm{ | } char(R\slash M) = p \textrm{ for some maximal ideal $M\subseteq R$}\}.  
\]
Now, for every subset $L$ of ${\mathbb P}\cup \{0\}$ there is a geometric formula $\phi_{L}$ in the empty context asserting that the characteristic (of a finitely presented von Neumann regular ring) is contained in $L$. In \cite{OC4} it is shown that, given $A$ and $B$ complementary subsets of ${\mathbb P}$, $\phi_{A}\stackrel{{\mathbb T }'}{\sim} \neg^{{\mathbb T}'} \phi_{B}$ and $\phi_{B}\stackrel{{\mathbb T }'}{\sim} \neg^{{\mathbb T}'} \phi_{A}$ where ${\mathbb T}'$ is the quotient of ${\mathbb T}$ given by the theory of fields. Moreover, it is shown, by working inside the universal model of ${\mathbb T}'$, that $\phi_{A}\vee \phi_{B}\stackrel{{\mathbb T }'}{\not \sim} \top$; hence, Theorem \ref{teofond} gives that the property of finitely presented von Neumann regular rings to have characteristic $\{0\}$ is not definable by a geometric formula in the signature of $\mathbb T$ (cfr. also \cite{OC4}).

{\bf Acknowledgements:} I am very grateful to my Ph.D. supervisor Peter Johnstone for useful comments on a draft of this paper.\\

\end{document}